\documentclass[11pt,twoside]{amsart}
\usepackage{}
\usepackage{amsthm}
\usepackage{amsmath}
\usepackage{amssymb}
\usepackage{amsfonts}
\usepackage{latexsym}
\usepackage{enumitem}
\usepackage{mathrsfs}
\usepackage[all]{xy} \SelectTips{eu}{}
\usepackage{hyperref}
\usepackage{eucal}
\usepackage{xcolor}

\newtheorem{intthm}{Theorem}[]

\newtheorem*{intque*}{Question}
\newtheorem*{intexa*}{Example}
\newcommand{\numberseries}{\bfseries}   %Fontseries used for numbering
\newlength{\thmtopspace}                %Space above theorem
\newlength{\thmbotspace}                %Space below theorem
\newlength{\thmheadspace}               %Space after theorem label
\newlength{\thmindent}                  %For indenting
\newtheoremstyle{bfupright head,slanted body}
{\thmtopspace}{\thmbotspace}
{\slshape}{\thmindent}{\bfseries}{.}{\thmheadspace}
{{\numberseries \thmnumber{#2\;}}\thmnote{#3}}

\newtheoremstyle{bfupright head,upright body}
{\thmtopspace}{\thmbotspace}
{\upshape}{\thmindent}{\bfseries}{.}{\thmheadspace}
{{\numberseries \thmnumber{#2\;}}\thmnote{#3}}

\newtheoremstyle{fixed bf head,slanted body}
{\thmtopspace}{\thmbotspace}{\slshape}
{\thmindent}{\bfseries}{.}{\thmheadspace}
{{\numberseries \thmnumber{#2\;}}\thmname{#1}\thmnote{ (#3)}}

\newtheoremstyle{fixed bf head,upright body}
{\thmtopspace}{\thmbotspace}{\upshape}
{\thmindent}{\bfseries}{.}{\thmheadspace}
{{\numberseries \thmnumber{#2\;}}\thmname{#1}\thmnote{ (#3)}}

\newtheoremstyle{numbered paragraph}
{\thmtopspace}{\thmbotspace}{\upshape}
{\thmindent}{\upshape}{}{\thmheadspace}
{{\numberseries \thmnumber{#2.}}}

\theoremstyle{bfupright head,slanted body}
\newtheorem{res}{}[section]             \newtheorem*{res*}{}

\theoremstyle{bfupright head,upright body}
\newtheorem{bfhpg}[res]{}               \newtheorem*{bfhpg*}{}

\theoremstyle{fixed bf head,slanted body}
\newtheorem{thm}[res]{Theorem}         \newtheorem*{thm*}{Theorem}
\newtheorem{prp}[res]{Proposition}      \newtheorem*{prp*}{Proposition}
\newtheorem{cor}[res]{Corollary}        \newtheorem*{cor*}{Corollary}
\newtheorem{lem}[res]{Lemma}            \newtheorem*{lem*}{Lemma}
         \newtheorem*{que*}{Question}
\theoremstyle{fixed bf head,upright body}
\newtheorem{dfn}[res]{Definition}       \newtheorem*{dfn*}{Definition}
\newtheorem{rmk}[res]{Remark}           \newtheorem*{rmk*}{Remark}
\newtheorem{exa}[res]{Example}           \newtheorem*{exa*}{Example}
\newtheorem{nota}[res]{Notation}           \newtheorem*{nota*}{Notation}
           \newtheorem*{setup*}{Setup}
\newtheorem{setup and notation}[res]{Setup and notation}  \newtheorem*{setup and notation*}{Setup and notation}
\theoremstyle{numbered paragraph}
\newtheorem{ipg}[res]{}

\newlength{\thmlistleft}        %leftmargin
\newlength{\thmlistright}       %rightmargin
\newlength{\thmlistpartopsep}   %partopsep
\newlength{\thmlisttopsep}      %topsep
\newlength{\thmlistparsep}      %parsep
\newlength{\thmlistitemsep}     %itemsep

\newcounter{eqc}
	{\end{list}}%

\newcounter{prt}
\newenvironment{prt}{\begin{list}{\upshape (\alph{prt})}%
		{\usecounter{prt}%
			\setlength{\leftmargin}{\thmlistleft}%
			\setlength{\labelwidth}{\thmlistleft}%
			\setlength{\rightmargin}{\thmlistright}%
			\setlength{\partopsep}{\thmlistpartopsep}%
			\setlength{\topsep}{\thmlisttopsep}%
			\setlength{\parsep}{\thmlistparsep}%
			\setlength{\itemsep}{\thmlistitemsep}}}%
	{\end{list}}%

\newcounter{rqm}
	{\end{list}}%

	{\end{list}}%

\newenvironment{prf*}[1][Proof]{%
	\begin{proof}[\bf #1]
		\setcounter{equation}{0}
		}
	{\end{proof}
}

\newcommand{\pgref}[1]{\ref{#1}}

\renewcommand{\eqref}[1]{(\pgref{eq:#1})}

\newcommand{\thmcite}[2][?]{\cite[Theorem~#1]{#2}}
\newcommand{\prpcite}[2][?]{\cite[Proposition~#1]{#2}}

\newcommand{\corcite}[2][?]{\cite[Corollary~#1]{#2}}
\newcommand{\lemcite}[2][?]{\cite[Lemma~#1]{#2}}

\newcommand{\dfncite}[2][?]{\cite[Definition~#1]{#2}}

\numberwithin{equation}{res}

\def\urltilda{\kern -.15em\lower .7ex\hbox{\~{}}\kern .04em}

\DeclareMathOperator*{\colim}{colim}

\newcommand{\GF}[1]{\mathsf{GF}(#1)}

\newcommand{\PGF}[1]{\mathsf{PGF}(#1)}
\newcommand{\GI}[1]{\mathsf{GI}(#1)}

\newcommand{\F}[1]{\mathsf{Flat}(#1)}

\newcommand{\AU}[1]{\mathsf{End}_{\calI}(#1)}

\newcommand{\Cot}[1]{\mathsf{Cot}(#1)}
\newcommand{\xra}[2][]{\xrightarrow[#1]{\:#2\:}}

\newcommand{\Prj}[1]{\mathsf{Proj}(#1)}
\newcommand{\Inj}[1]{\mathsf{Inj}(#1)}

\newcommand{\Hom}[3][]{\operatorname{Hom}_{#1}(#2,#3)}

\newcommand{\tp}[3][{\sf E}]{\nobreak{#3\otimes_{#1}#2}}
\newcommand{\Ext}[4][R]{\operatorname{Ext}_{#1}^{#2}(#3,#4)}

\newcommand{\POMon}{\mathsf{POMon}}
\newcommand{\PBEpi}{\mathsf{PBEpi}}
\newcommand{\Mon}{\mathsf{Mon}}
\newcommand{\Epi}{\mathsf{Epi}}

\newcommand{\eva}{\mathsf{eva}}
\newcommand{\fre}{\mathsf{fre}}
\newcommand{\cofre}{\mathsf{cofre}}

\newcommand{\calA}{\mathcal{A}}

\newcommand{\calC}{\mathcal{C}}
\newcommand{\calD}{\mathcal{D}}
\newcommand{\calF}{\mathcal{F}}
\newcommand{\calI}{\mathcal{I}}
\newcommand{\calJ}{\mathcal{J}}

\newcommand{\calP}{\mathcal{P}}
\newcommand{\calQ}{\mathcal{Q}}
\newcommand{\calR}{\mathcal{R}}
\newcommand{\calS}{\mathcal{S}}

\newcommand{\calV}{\mathcal{V}}
\newcommand{\calW}{\mathcal{W}}
\newcommand{\calX}{\mathcal{X}}
\newcommand{\calY}{\mathcal{Y}}
\newcommand{\sfI}{\mathsf{I}}
\newcommand{\sfJ}{\mathsf{J}}
\newcommand{\sfC}{\mathsf{C}}

\newcommand{\sfF}{\mathsf{F}}
\newcommand{\sfL}{\mathsf{L}}
\newcommand{\sfM}{\mathsf{M}}
\newcommand{\sfW}{\mathsf{W}}
\newcommand{\bbZ}{\mathbb{Z}}
\newcommand{\bbQ}{\mathbb{Q}}
\newcommand{\scrD}{\mathscr{D}}
\newcommand{\scrR}{\mathscr{R}}

\newcommand{\Mor}{\mathsf{Mor}}
\newcommand{\Ob}{\mathsf{Ob}}
\newcommand{\lMod}{\textrm{-} \mathsf{Mod}}
\newcommand{\rMod}{\mathsf{Mod} \textrm{-}}
\newcommand{\lRep}{\textrm{-} \mathsf{Rep}}

\newcommand{\K}{\mbox{\sf ker}}
\newcommand{\coker}{\mbox{\rm coker}}

\newcommand{\op}{^{\sf op}}
\newcommand{\qis}{\simeq}
\newcommand{\C}{\mbox{\sf cok}}
\newcommand{\id}{\mathrm{id}}

\newcommand{\qra}{\xra{\qis}}

\def\soft#1{\leavevmode\setbox0=\hbox{h}\dimen7=\ht0\advance
	\dimen7 by-1ex\relax\if t#1\relax\rlap{\raise.6\dimen7
		\hbox{\kern.3ex\char'47}}#1\relax\else\if T#1\relax
	\rlap{\raise.5\dimen7\hbox{\kern1.3ex\char'47}}#1\relax
	\else\if d#1\relax\rlap{\raise.5\dimen7\hbox{\kern.9ex
			\char'47}}#1\relax\else\if D#1\relax\rlap{\raise.5\dimen7
		\hbox{\kern1.4ex\char'47}}#1\relax\else\if l#1\relax
	\rlap{\raise.5\dimen7\hbox{\kern.4ex\char'47}}#1\relax
	\else\if L#1\relax\rlap{\raise.5\dimen7\hbox{\kern.7ex
			\char'47}}#1\relax\else\message{accent \string\soft
		\space #1 not defined!}#1\relax\fi\fi\fi\fi\fi\fi}
\makeatletter
\def\part{\@startsection{part}{1}%
\z@{.7\linespacing\@plus\linespacing}{.8\linespacing}%
{\LARGE\sffamily\centering}}
\makeatletter
\def\l@section{\@tocline{1}{2pt}{0pc}{}{}}
\makeatother
\let\oldtocpart=\tocpart
\renewcommand{\tocpart}[2]{\bf\large\oldtocpart{#1}{#2}}
\let\oldtocsection=\tocsection
\renewcommand{\tocsection}[2]{\bf\oldtocsection{#1}{#2}}

\title[Representations over diagrams of abelian categories]{Representations over diagrams of abelian categories II: Abelian model structures}



\author[Z.X. Di]{Zhenxing Di}
\address{Z.X. Di \ School of Mathematical Sciences, Huaqiao University, Quanzhou 362021, China}
\email{dizhenxing@163.com}

\author[L.P. Li]{Liping Li}
\address{L.P. Li \ LCSM (Ministry of Education), Department of Mathematics, Hunan Normal University, Changsha 410081, China}
\email{lipingli@hunnu.edu.cn}

\author[L. Liang]{Li Liang}
\address{L. Liang \
Department of Mathematics, Gansu Center for Fundamental Research in Complex Systems Analysis and Control, Lanzhou Jiaotong University, Lanzhou 730070, China}
\email{lliangnju@gmail.com}
\urladdr{https://sites.google.com/site/lliangnju}

\author[N.N. Yu]{Nina Yu}
\address{N.N. Yu \ School of Mathematical Sciences, Xiamen University, Xiamen 361005, China}
\email{ninayu@xmu.edu.cn}

\begin{document}

\begin{abstract}
This is the second paper in a series on representations over diagrams of abelian categories. We show that, under certain conditions, a compatible family of abelian model categories indexed by a skeletal small category can be amalgamated into an abelian model structure on the category of representations. Our approach focuses on classes of morphisms rather than cotorsion pairs of objects. Additionally,  we provide an explicit description of cofibrant objects in the resulting abelian model category. As applications, we construct Gorenstein injective and Gorenstein flat model structures on the category of presheaves of modules over a special class of index category and characterize  Gorenstein homological objects within this framework.
\end{abstract}

%%%%%%%%%%%%%%%%%%%%%%%%%%%%%%%%%%%%%%%
\maketitle
\vspace*{.3cm}
\tableofcontents
\vspace*{-.9cm}
\enlargethispage{.4cm}
%%%%%%%%%%%%%%%%%%%%%%%%%%%%%%%%%%%%%%

\section*{Introduction}
\noindent
Throughout the paper, we let $\calI$ be a skeletal small category with the set of objects $\Ob(\calI)$ and the set of morphisms $\Mor(\calI)$. In the first paper of this series  \cite{DLLY}, we studied diagrams $\scrD$ of abelian categories defined as pseudo-functors from $\calI$ (viewed as a 2-category) to the meta-2-category of abelian categories, and representations $M$ over $\scrD$ assigning to each object $i$ in $\calI$ an object $M_i$ in the abelian category $\scrD_i$ such that certain compatibility conditions are satisfied. We refer the read to \cite{DLLY} for further details on  diagrams of abelian categories and their representations.
In particular, we investigated the Grothendieck structure of the category $\scrD \lRep$ of representations over $\scrD$, and characterized special homological objects such as projective  and injective objects. The main goal of the present paper is to construct various abelian model structures on $\scrD \lRep$ and to give explicit descriptions of special classes of objects, including cofibrant ones.

In this paper, we continue to use  $\scrD$ to denote  an $\calI$-diagram of Grothendieck categories admitting enough projectives, such that the functor $\scrD_{\alpha}: \scrD_i \to \scrD_j$ is right exact and preserves small coproducts for every morphism $\alpha: i \to j$ in $\calI$. Under these assumptions, the category $\scrD \lRep$ is a Grothendieck category admitting enough projectives, as shown in  \cite[Theorem 2.8]{DLLY}.

Model category theory, introduced by Quillen in \cite{Qui67}, provides a foundational framework for incorporating homotopy theory into  categorical settings. %Given a bicomplete category $\calA$ and a class $\calW$ of morphisms in $\calA$, one may wish to treat the morphisms in $\calW$ as isomorphisms, that is, to understand the localization category $\calA[\calW^{-1}]$. By imposing a certain model structure on $\calA$ one can construct a new category ${\rm Ho}(\calA)$, the homotopy category of $\calA$. A fundamental result shows that there is a canonical equivalence $\calA[\calW^{-1}] \cong {\rm Ho}(\calA)$; see \cite{Ho99}.
A model category is called \emph{abelian} \cite{Ho02} if its underlying category $\calA$ is abelian and the model structure is compatible with the abelian structure of $\calA$.
A celebrated result bridging  complete cotorsion pairs and abelian model structures was established by Hovey in \cite{Ho02}, now known as \emph{Hovey's correspondence}. It asserts that an abelian model structure on a bicomplete abelian category $\calA$ corresponds bijectively to a triple $(\calQ, \calW, \calR)$ of subcategories of $\calA$, where  $\calW$ is a thick subcategory, and both $(\calQ, \calW \cap \calR)$ and $(\calQ \cap \calW, \calR)$ are complete cotorsion pairs  in $\calA$. Here $\calQ$,  $\calW$ and $\calR$ respectively consist of cofibrant, trivial  and fibrant objects associated to the corresponding abelian model structure. Hovey's correspondence thus allows  an abelian model structure on $\calA$ to be succinctly represented by the triple $(\calQ, \calW, \calR)$, which is commonly referred to in the literature as a \emph{Hovey triple}. For more details on abelian model structures, we refer the reader to  \cite{Be14, Gil16, Gil162, Ho02}.

Suppose that each category $\scrD_i$ admits an abelian model structure, and that these  structures are compatible with respect to the diagram $\scrD$. A natural question arises:  Can these local model structures be amalgamated into an abelian model structure on the category of representations $\scrD \lRep$.  In the special case where $\scrD$ is a trivial diagram, this question reduces to whether an abelian model structure on $\calA$  induces one on the functor category $\mathrm{Fun}(\calI,\calA)$. This problem has been  investigated by several authors using Hovey's correspondence; see for example,  \cite{DELO,HJ19b}. However, a major obstacle in this approach is the difficulty of  verifying  the completeness of the induced cotorsion pairs in $\mathrm{Fun}(\calI,\calA)$.

We approach this question from a different perspective. For a special type of index category, the works in \cite{Ho99} and \cite{HR08} provide a method for constructing a model structure on $\scrD \lRep$ from a compatible family of model categories $\scrD_i$, by  focusing on morphisms rather than objects. This naturally leads to the question of  whether the resulting model structure  on $\scrD \lRep$ is abelian, assuming  that each  $\scrD_i$ carries an  abelian model structure. A key  advantage of this approach is that it does not rely on the Hovey's correspondence, thereby avoiding the challenge of verifying the completeness of induced cotorsion pairs.

As the first main result of this approach, we construct cofibrantly generated (hereditary) Hovey triples in $\scrD \lRep$ from cofibrantly generated (hereditary) Hovey triples in $\scrD_i$; see Theorem \ref{ht induce ht 2}.

\begin{intthm}\label{thmC}
Suppose that $\scrD$ is exact, and that $(\calQ_i, \calW_i, \calR_i)$ is a cofibrantly generated (hereditary) Hovey triple in $\scrD_i$ for each  $i \in \Ob(\calI)$. Assume further that both $\{ \calQ_i \}_{i \in \Ob(\calI)}$ and $\{\calQ_i \cap \calW_i\}_{i \in \Ob(\calI)}$ are compatible with respect to $\scrD$. Then
\[
(^\perp(\scrD \lRep_{\calW \cap \calR}), \, \scrD \lRep_{\calW}, \, \scrD \lRep_{\calR})
\]
is a cofibrantly generated (hereditary) Hovey triple in $\scrD \lRep$, where $\scrD \lRep_{\calW}$ (resp., $\scrD \lRep_{\calR}$ and $\scrD \lRep_{\calW \cap \calR}$) denotes the full subcategory of $\scrD \lRep$ consisting of representations $M$ such that $M_i\in \calW_i $ (resp., $\calR_i$ and $\calW_i \cap \calR_i$) for all $i \in \Ob(\calI)$.
\end{intthm}

A limitation of Theorem~\ref{thmC} is that the class $^\perp(\scrD \lRep_{\calW \cap \calR})$ is generally difficult to describe explicitly, even though the result applies to arbitrary index categories $\calI$. For a special class of index categories known as  \textit{direct categories}, H\"{u}ttemann and R\"{o}ndigs \thmcite[3.2.13]{HR08} constructed a model structure on $\scrD \lRep$, induced from a compatible family of model structures on the categories $\scrD_i$, and provided an explicit description of the cofibrations. This naturlly raises the question of whether their construction extends to the setting of  abelian model structures. Unfortunately, it appears that this is not ture  in general. However, under certain conditions, we are able to obtain an affirmative answer.

Inspired by the works of Auslander and Reiten \cite{AR91} and Birkhoff \cite{B34},  for a family $\calX = \{ \calX_i \}_{i \in \Ob(\calI)}$ of subcategories of $\scrD_i$, we define a subcategory $\Phi(\calX)$ of $\scrD \lRep$ (see Definition \ref{s and phi}) whose objects admit a transparent description. In particular,  $\Phi(\calX)$ coincides with the monomorphism category introduced in \cite{Z11}. Under suitable conditions, we show that $$^\perp(\scrD \lRep_{\calW \cap \calR}) = \Phi(\calQ),$$  leading to the following result; see Theorem \ref{ht induce ht 1}.

\begin{intthm}\label{thmB}
Suppose that $\calI$ is a left rooted quiver (viewed as a category in a natural way), and that $\scrD$ is exact. If $(\calQ_i, \calW_i, \calR_i)$ is a (hereditary) Hovey triple in $\scrD_i$ for each $i\in\Ob(\calI)$, and the families
 $\calQ = \{ \calQ_i \}_{i \in \Ob(\calI)}$ and $\{\calQ_i \cap \calW_i\}_{i \in \Ob(\calI)}$ are compatible with respect to $\scrD$, then there exists a (hereditary) Hovey triple in $\scrD \lRep$ given by
\[
(\Phi(\calQ),\, \scrD \lRep_{\calW},\, \scrD \lRep_{\calR}).
\]
\end{intthm}

We then give some applications of Theorem \ref{thmB}. For finitely generated modules over a noetherian ring, Auslander and Bridger \cite{AusBri} introduced the \textit{G-dimension}, a homological invariant that was later generalized to arbitrary modules. In \cite{EEnOJn95b, EJT-93}, Enochs, Jenda and Torrecillas introduced two primary generalized modules of G-dimension zero: Gorenstein projective modules and Gorenstein flat modules.  These notions have proven to be particularly important and were further developed by Holm in \cite{HHl04a}. For an  associative ring $A$, \v{S}aroch and \v{S}t'ov\'{\i}\v{c}ek  \cite{SS20} introduced projectively coresolved Gorenstein flat $A$-modules and use them to construct a new abelian model structure on $A \lMod$, the category of left $A$-modules. Specifically, let $\GF{A}$ and  $\PGF{A}$ denote the subcategory of Gorenstein flat and projectively coresolved Gorenstein flat $A$-modules, respectively. They show that there exists a hereditary Hovey triple $(\GF{A}, \PGF{A}^\perp, \Cot{A})$ in $A \lMod $, where $\Cot{A}$ is the subcategory of cotorsion $A$-modules.

Let $\scrR$ be a representation of $\calI$ on the category {\sf Ring}, assigning each $i\in\calI$ an associative ring $\scrR_i=R_i$ and to each $\alpha \in \Mor(\calI)$  a ring homomorphism $\scrR_\alpha: R_i\to R_j$. Then \v{S}aroch and \v{S}t'ov\'{\i}\v{c}ek' s result yields a hereditary Hovey triple
$(\GF{R_i}, \, \PGF{R_i}^\perp, \, \Cot{R_i})$
in $R_i \lMod $ for every $i \in \Ob(\calI)$. We mention that the category of left $\scrR$-modules studied by Estrada and Virili in \cite{SS2017} coincides with $\overline{\scrR} \lRep$, where $\overline{\scrR}$ is a right exact $\calI$-diagram of left module categories induced by $\scrR$ (see Example \ref{bar R}).
By Lemma \ref{GF and PGF compatible}, the families ${\sf Flat_{\bullet}} = \{ \F{R_i} \}_{i \in \Ob(\calI)}$ and ${\sf GF_{\bullet}} = \{ \GF{R_i} \}_{i \in \Ob(\calI)}$  are compatible with respect to $\overline{\scrR}$. Thus as an application of Theorem \ref{thmB}, we get a hereditary Hovey triple in the category $\overline{\scrR} \lRep$ and classify all Gorenstein flat objects and cotorsion objects in $\overline{\scrR} \lRep$; for details, see Theorems \ref{ht induce ht 12} and \ref{Gorenstein flat r module}, Proposition \ref{cotorsion object} and Corollary \ref{cp for pgf}.

\begin{intthm} \label{thmG}
Let $\calI$ be a left rooted quiver and $\scrR$ a flat representation of $\calI$ on {\sf Ring}.
Then there exists a hereditary Hovey triple
\[
(\Phi({\sf GF_{\bullet}}),\ \overline{\scrR} \lRep_{{\sf PGF}^\perp_{\bullet}},\ \overline{\scrR} \lRep_{\sf Cot_{\bullet}})
\]
in $\overline{\scrR} \lRep$, where:
\begin{itemize}
    \item
$\Phi({\sf GF_{\bullet}})$ coincides with the subcategory consisting of Gorenstein flat objects in $\overline{\scrR} \lRep$, \item  $\overline{\scrR} \lRep_{{\sf PGF}^\perp_{\bullet}}$ coincides with the right orthogonal subcategory $\PGF{\overline{\scrR} \lRep}^\perp$,
\item  $\overline{\scrR} \lRep_{\sf Cot_{\bullet}}$ coincides with the subcategory consisting of cotorsion objects in $\overline{\scrR} \lRep$.
\end{itemize}
\end{intthm}

We can also construct a Gorenstein injective model structure on $\widetilde{\scrR} \lRep$ and provide an explicit description of the  Gorenstein injective objects in this category (see Theorems \ref{ht induce ht GI 1} and \ref{Gorenstein in r module} and Corollary \ref{cp for ginj}). Here,  $\widetilde{\scrR}$ denotes the  exact $\calI\op$-diagram of right module categories induced by the representation $\scrR$ of $\calI$ on {\sf Ring} (see Example \ref{tilder R}).

\begin{intthm} \label{thmD}
Let $\calI$ be a left rooted quiver and $\scrR$ a flat representation of $\calI$ on {\sf Ring}. Then there exists a hereditary Hovey triple
\[
(\widetilde{\scrR} \lRep, \, \widetilde{\scrR} \lRep_{^\perp \sf GI_{\bullet}}, \, \Psi(\sf GI_{\bullet}))
\]
in $\widetilde{\scrR} \lRep$, where
\begin{itemize}
    \item
$\Psi(\sf GI_{\bullet})$ coincides with the subcategory consisting of Gorenstein injective objects in $\widetilde{\scrR} \lRep$,
  \item  $\widetilde{\scrR} \lRep_{^\perp \sf GI_{\bullet}}$ coincides with the right orthogonal subcategory ${^\perp\GI{\widetilde{\scrR} \lRep}}$.
\end{itemize}
\end{intthm}

\section{Abelian model structures}
\label{Preliminary on abelian model structures}
\noindent
We begin this section by reviewing some  background  on abelian model category theory.

\begin{bfhpg}[\bf Weak factorization system]
Let $l: A \to B$ and $r: C \to D$ be morphisms in a category $\calD$. Recall that $l$ has the {\it left lifting property} with respect to $r$ (or $r$ has the {\it right lifting property} with respect to $l$) if for every pair of morphisms $f: A \to C$ and $g: B \to D$ with $rf = gl$, there exists a morphism $t: B\to C$ such that the diagram
\[
\xymatrix{
A \ar[d]_{l} \ar[r]^{f} & C \ar[d]^{r}\\
B \ar[ur]|-{t} \ar[r]_{g} & D
}
\]
commutes.

For a class $\sfC$ of morphisms in $\calD$, let $\sfC^\Box$ denote the class of morphisms in $\calD$ that have the right lifting property with respect to all morphisms in $\sfC$. The class $^\Box\sfC$ is defined dually. Recall from Bousfield \cite{Bo77} that a pair $(\sfC, \sfF)$ of classes of morphisms in $\calD$ is called a {\it weak factorization system} if \begin{itemize}\item $\sfC^\Box = \sfF$ and ${^\Box\sfF} = \sfC$, and
\item every morphism $\alpha$ in $\calD$ can be factored as $\alpha = fc$ with $c \in \sfC$ and $f \in \sfF$.
\end{itemize}
\end{bfhpg}

The following definition of model categories is slightly more general than that given by Hovey \cite{Ho99}, in which the factorizations are required to be functorial.

\begin{dfn}
Let $\calD$ be a bicomplete category. A {\it model structure} on $\calD$ is a triple $(\sfC, \sfW, \sfF)$ of classes of morphisms in $\calD$ such that:
\begin{prt}
\item $(\sfC, \sfW \cap \sfF)$ and $(\sfC \cap \sfW, \sfF)$ are weak factorization systems;

\item $\sfW$ satisfies the 2-out-of-3 property: if two of the three morphisms $\alpha$, $\beta$ and $\beta \alpha$ lie in $\sfW$, then so does the third.
\end{prt}

Morphisms in $\sfC$ (resp., $\sfW$, $\sfF$) are called {\it cofibrations} (resp., {\it weak equivalences}, {\it fibrations}). Morphisms in $\sfW \cap \sfF$ (resp., $\sfC \cap \sfW$) are called {\it trivial fibrations} (resp., {\it trivial cofibrations}).
A model structure $(\sfC, \sfW, \sfF)$ on $\calD$ is called {\it cofibrantly generated} if there exist sets $\sfI$ (of generating cofibrations) and $\sfJ$ (of generating trivial cofibrations) of morphisms  such that both $I$ and $J$ permit the small object argument and $\sfI^\Box=\sfW\cap\sfF$ and $\sfJ^\Box=\sfF$. An object in $\calD$ is called \emph{cofibrant} if the morphism from the initial object to it  is a cofibration. Fibrant objects are defined dually. An object in $\calD$ is called \emph{trivial} if the morphism from the initial object to it  is a weak equivalence, or equivalently, the morphism from it to the terminal object is a weak equivalence.
\end{dfn}

\begin{dfn}\label{df of abelian ms}
Let $\calA$ be a bicomplete abelian category. Recall from \cite{Ho02} that a model structure on $\calA$ is said to be \emph{abelian} if the following conditions are satisfied:
\begin{prt}
\item Every cofibration is a monomorphism;

\item Every fibration is an epimorphism with fibrant kernel;

\item Every Trivial fibration is an epimorphism with trivial fibrant kernel.
\end{prt}
\end{dfn}

\begin{rmk}\label{df of ams is self-dual}
Let $\calA$ be a bicomplete abelian category. By \prpcite[4.2]{Ho02}, a model structure on $\calA$ in which cofibrations are monomorphisms and fibrations are epimorphisms is abelian if and only if:
\begin{itemize} \item Cofibrations coincide with monomorphisms whose cokernels are cofibrant;
\item Trivial cofibrations coincide with monomorphisms whose cokernels are trivial and cofibrant.
\end{itemize}
\end{rmk}

\begin{bfhpg}[\bf Cotorsion pairs]\label{cotpair}
The concept of cotorsion pairs was first introduced by Salce \cite{S-L} and rediscovered by Enochs and Jenda in \cite{rha}. It  is an analogue of  torsion pairs, with the Hom functor replaced by the Ext functor. Let $\calA$ be an abelian category. A pair $(\calC,\calF)$ of subcategories of $\calA$ is called a \emph{cotorsion pair} if $$\calC^{\bot} = \calF\quad \text{and}\quad  ^{\bot}\calF =\calC,$$ where
\begin{itemize}
\item $\calC^{\bot} = \{M \in\calA \mid \textrm{Ext}_{\calA}^{1}(C, M) = 0 \textrm{ for all objects } C \in \calC \}$;

\item $^{\bot}\calF = \{M \in \calA \mid \textrm{Ext}_{\calA}^{1}(M, D) = 0
       \textrm{ for all objects } D \in \calF\}$.
\end{itemize}

Following \cite{rha}, a cotorsion pair $(\calC,\calF)$ is said to be \emph{cogenerated by a set} if there is a set $\calS$ of objects in $\calA$ such that $\calS^{\bot}=\calF$. A cotorsion pair $(\calC,\calF)$ is called \emph{complete} if for every object $M$ in $\calA$, there exist short exact sequences $0 \to D \to C \to M \to 0$ and $0 \to M \to D' \to C' \to 0$ in $\calA$ with $D, D' \in \calF$ and $C, C' \in \calC$. A cotorsion pair $(\calC,\calF)$ is called \textit{resolving} if $\calC$ is closed under taking kernels of epimorphisms between objects of $\calC$, and \textit{coresolving} if $\calF$ is closed under taking cokernels of monomorphisms between objects in $\calF$. A cotorsion pair is called \emph{hereditary} if it is both resolving and coresolving.
\end{bfhpg}

The following result is due to Becker \cite{Be14}.

\begin{lem}\label{test hereditary}
Let $\calA$ be an abelian category and $(\calC,\calF)$ a complete cotorsion pair in $\calA$. Then $(\calC,\calF)$ is hereditary if and only if it is resolving, if and only if it is coresolving.
\end{lem}

The central result in the theory of abelian model categories is now known as  Hovey's correspondence. Recall that a subcategory of an abelian category $\calA$ is called {\it thick} if it is closed under direct summands, extensions, kernels of epimorphisms and cokernels of monomorphisms. In what follows, for a subcategory $\calY$ of $\calA$, set
\begin{itemize}
\item $\Mon(\calY)=\{ f \ |\ f \text{ is a monomorphism with } \coker(f) \in \calY\}$
\item $\Epi(\calY)=\{ f \ |\ f \text{ is an epimorphism with } \ker(f) \in \calY\}$.
\end{itemize}

\begin{thm}[Hovey's correspondence]\label{Hovey cor}
Let $\calA$ be a bicomplete abelian category. Then there exists a bijective correspondence between
\begin{prt}
\item abelian model structures $(\sfC, \sfW, \sfF)$ on $\calA$, and

\item triples $(\calQ, \calW, \calR)$ of subcategories of $\calA$ such that both $(\calQ, \calW \cap \calR)$ and $(\calQ \cap \calW, \calR)$ are complete cotorsion pairs in $\calA$, and $\calW$ is thick.
\end{prt}
Explicitly, given an abelian model structure $(\sfC, \sfW, \sfF)$ on $\calA$, the corresponding triple of subcategories of $\calA$ consists of the cofibrant, trivial, and fibrant objects, respectively. Conversely, given a triple $(\calQ, \calW, \calR)$ of subcategories of $\calA$ as in {\rm (b)}, the associated abelian model structure is $(\Mon(\calQ), \sfW, \Epi(\calR))$, where
\[
\sfW = \{ w \,|\, w \text{ can be decomposed as } w = fc \text{ with } c \in \Mon(\calQ \cap \calW) \text{ and } f \in \Epi(\calW \cap \calR) \}.
\]
\end{thm}

Hovey's correspondence shows that an abelian model structure on $\calA$ can be succinctly represented by a triple of subcategories of $\calA$ satisfying the conditions in (b). Such a triple is often referred to  as an \emph{abelian model structure}, and is called a \emph{Hovey triple}.

\begin{dfn}
A Hovey triple $(\calQ, \calW, \calR)$ in a bicomplete abelian category $\calA$ is said to be \emph{cofibrantly generated} if the associated abelian model structure $(\Mon(\calQ), \sfW, \Epi(\calR))$ on $\calA$ is cofibrantly generated. It is said to be  \emph{hereditary} if both the cotorsion pairs $(\calQ, \calW \cap \calR)$ and $(\calQ \cap \calW, \calR)$  are hereditary.
\end{dfn}

\begin{rmk}
Let $(\calQ, \calW, \calR)$ be a Hovey triple in a Grothendieck category $\calA$ with enough projectives. If both cotorsion pairs $(\calQ, \calW \cap \calR)$ and $(\calQ \cap \calW, \calR)$  are cogenerated by sets, then it follows from \cite[Lemma 6.7 and Corollary 6.8]{Ho02} that the Hovey triple $(\calQ, \calW, \calR)$ is cofibrantly generated. Indeed, most of Hovey triples arising in practice are cofibrantly generated.
\end{rmk}

The thick subcategory $\calW$ in a Hovey triple $(\calQ, \calW, \calR)$ plays a central role,  as
it determines the homotopy category of the associated abelian model category.
The following result provides a characterization of $\calW$; see Gillespie \prpcite[3.2]{Gil16}.

\begin{lem} \label{model structor111}
Let $(\calQ,\calW,\calR)$ be a Hovey triple in a bicomplete abelian category $\calA$. Then the thick subcategory $\calW$ can be characterized as
\begin{align*}
\calW & = \{M \mid \text{ there is a s.e.s. } 0 \to M \to A \to B \to 0 \text{ with } A \in \calW \cap \calR \text{ and } B \in \calQ \cap \calW \}\\
& = \{ M \mid \text{ there is a s.e.s. } 0 \to A' \to B' \to M \to 0 \text{ with } A' \in \calW \cap \calR \text{ and } B' \in \calQ \cap \calW \}.
\end{align*}
Consequently, $\calW$ is uniquely determined by $\calQ$ and $\calR$; that is, if  $(\calQ,\mathcal{V},\calR)$ is another Hovey triple, then $\mathcal{V} = \calW$.
\end{lem}

Additional background  on abelian model structures can be found in \cite{Gil162}.

\section{Representations over diagrams of abelian categories}
\noindent
In this section, we provide background  on representations over diagrams of categories. For further details, the read may  refer to \cite{DLLY}.

\begin{dfn} \label{generalised diagram}
An $\calI$-\emph{diagram} of categories is a tuple $(\scrD,\eta,\tau)$ (often simply denoted by $\scrD$) consisting of the following data:
\begin{itemize}
\item For each $i \in \Ob(\calI)$, a category $\scrD_i$;
\item For each $\alpha: i \to j \in \Mor(\calI)$, a covariant functor $\scrD_\alpha: \scrD_i \to \scrD_j$;
\item For each $i \in \Ob(\calI)$, a natural isomorphism $\eta_i: {\id}_{\scrD_i} \qra \scrD_{{e}_i}$, where ${e}_i$ is the identity on $i$;
\item For each pair of composable morphisms $\alpha$ and $\beta$ in $\Mor(\calI)$, a natural isomorphism
\[
\tau_{\beta,\alpha}: \scrD_{\beta} \circ \scrD_{\alpha} \qra \scrD_{\beta\alpha}
\]
\end{itemize}
such that the following two axioms are satisfied:

\noindent ({Dia.1}) Given composable morphisms $i \overset{\alpha} \to j \overset{\beta} \to k \overset{\gamma} \to l \in \Mor(\calI)$, the diagram
\[
\xymatrix{
\scrD_{\gamma} \circ \scrD_{\beta} \circ \scrD_{\alpha} \ar[rr]^-{\id_{\scrD_\gamma} \ast \tau_{\beta,\alpha}} \ar[d]_{\tau_{\gamma,\beta} \ast \id_{\scrD_\alpha}} & & \scrD_{\gamma} \circ \scrD_{\beta\alpha} \ar[d]^{\tau_{\gamma, \beta\alpha}}\\
\scrD_{\gamma\beta} \circ \scrD_{\alpha} \ar[rr]^-{\tau_{\gamma\beta, \alpha}} & & \scrD_{\gamma\beta\alpha}
}
\]
of natural isomorphisms commutes, where ``$\ast$'' denotes the Godement product of natural transformations.

\noindent ({Dia.2}) Given a morphism $i \overset{\alpha} \to j \in \Mor(\calI)$, the diagram
\[
\xymatrix{
 & \scrD_{\alpha} \ar@{=}[dd] \ar[dl]_{\id_{\scrD_{\alpha}} \ast \eta_i} \ar[dr]^{\eta_j \ast \id_{\scrD_{\alpha}}} & & \\
\scrD_{\alpha} \circ \scrD_{e_i} \ar[dr]_{\tau_{\alpha,e_i}} & & \scrD_{e_j} \circ \scrD_{\alpha} \ar[dl]^{\tau_{e_j, \alpha}}\\
& \scrD_{\alpha}
}
\]
of natural isomorphisms commutes.
\end{dfn}

An $\calI$-diagram $\scrD$ of categories is said to be \emph{strict} if $\eta_i$ is the identity for any $i \in \Ob(\calI)$ and $\tau_{\beta,\alpha}$ is the identity for any pair of composable morphisms $\alpha$ and $\beta$ in $\Mor(\calI)$. It  is said to be \emph{admitting enough right adjoints} if each functor $\scrD_\alpha$ admits a right adjoint for any $\alpha \in \Mor(\calI)$; it is called an adjunction bundle or $\calI$-bundle by H\"{u}ttemann and R\"{ondigs} in \cite{HR08}.

An $\calI$-diagram $\scrD$ of abelian categories is called \emph{exact} (resp., \emph{right exact}) if each functor $\scrD_\alpha$ is exact (resp., right exact) and additive. For example, let $\calA$ be an abelian category. Then the diagram defined by  $\scrD_i = \calA$ for all $i \in \Ob(\calI)$ and $\scrD_\alpha = \id_{\calA}$ for all $\alpha\in\Mor(\calI)$ is an  exact $\calI$-diagram, called the \emph{trivial} $\calI$-diagram of $\calA$.

\begin{dfn} \label{DF OF R-M}
Let $(\scrD,\eta,\tau)$ be an $\calI$-diagram of categories. A \textit{representation} $M$ over $\scrD$ consists of the following data:
\begin{itemize}
\item For each $i \in \Ob(\calI)$, an object $M_i \in \scrD_i$;
\item For each morphism $\alpha: i \to j \in \Mor(\calI)$,
a structural morphism $M_\alpha: \scrD_{\alpha}(M_i) \to M_j \in \scrD_j$
\end{itemize}
such that the following two axioms hold:

\noindent (Rep.1) For any composable morphisms
$i \overset{\alpha} \to j \overset{\beta} \to k \in \Mor(\calI)$,
the diagram
\[
\xymatrix{
\scrD_{\beta\alpha}(M_i) \ar[rr]^-{M_{\beta\alpha}} & & M_k\\
\scrD_{\beta}(\scrD_{\alpha}(M_i)) \ar[u]^{\tau_{\beta,\alpha}(M_i)} \ar[rr]^-{\scrD_{\beta}({M_\alpha})} & & \scrD_{\beta}(M_j) \ar[u]_{M_\beta}
}
\]
in $\scrD_k$ commutes; that is, $M_{\beta\alpha} \circ \tau_{\beta,\alpha}(M_i) = {M_\beta} \circ \scrD_{\beta}({M_\alpha})$.

\noindent (Rep.2) For each $i \in \Ob(\calI)$, the diagram
\[
\xymatrix{
M_i \ar[rr]^-{\id_{M_i}} \ar[dr]_{\eta_i(M_i)} &  &  M_i\\
 & \scrD_{e_i}(M_i) \ar[ur]_{M_{e_i}}              }
\]
in $\scrD_i$ commutes;  that is, $M_{e_i} = \eta^{-1}_i(M_i)$.

A morphism $\omega: M \to M'$ between two representations $M$ and $M'$ over $\scrD$ is a family  of morphisms $\{ \omega_i: M_i \to M'_i\}_{i \in \Ob(\calI)}$ such that for every $\alpha: i \to j \in \Mor(\calI)$, the diagram
\[
\xymatrix{
\scrD_{\alpha}(M_i) \ar[rr]^-{\scrD_{\alpha}(\omega_i)} \ar[d]_{M_{\alpha}} & & \scrD_{\alpha}(M'_i) \ar[d]^{M'_{\alpha}}\\
M_j \ar[rr]^-{\omega_{j}} & & M'_j
}
\]
in $\scrD_j$ commutes.
\end{dfn}

Representations over $\scrD$ were referred to as twisted diagrams in \cite{HR08}. We denote by $\scrD \lRep$ the category of all representations over $\scrD$. According to \prpcite[1.8]{DLLY}, if $\scrD$ is a right exact $\calI$-diagram of abelian categories, then  $\scrD \lRep$ is an abelian category.  A sequence $M \to N \to K$ in $\scrD \lRep$ is exact if and only if $M_i \to N_i \to K_i$ is exact in $\scrD_i$ for each $i \in \Ob(\calI)$.
The category $\scrD \lRep$ provides a unifying framework for a wide range of categories, including  comma categories, module categories over Morita context rings, categories of additive functors from $\calI$ to an abelian category, and  categories of representations of (generalized) species and phyla; see \cite{DLLY}. In the following,  we give two examples that will be used in the final three sections of this paper.

\begin{exa}\label{bar R}
Estrada and Virili introduced in \cite{SS2017} the notion of representations $\scrR$ of $\calI$ on {\sf Ring}, the category of associative rings. These are, in our terminology,  $\calI$-diagrams of associative rings, where for each $i\in\calI$, the component $\scrR_i=R_i$ an associative ring (viewed as a preadditive small category with one object), and for each $\alpha \in \Mor(\calI)$,   $\scrR_\alpha$ is a ring homomorphism.  In \cite{SS2017}, the authors also introduced the notion of  modules over such representations and established several important homological properties of their module category. One can associate to a representation $\scrR$ of $\calI$ on {\sf Ring}  an $\calI$-diagram $\overline{\scrR}$ of left module categories defined by:
\begin{itemize}
\item $\overline{\scrR}_i = R_i \lMod$, the category of left $R_i$-modules,
      for each $i \in \Ob(\calI)$;
\item $\overline{\scrR}_\alpha = R_j \otimes_{R_i}-: R_i\lMod \to R_j\lMod$
      for each morphism $\alpha: i \to j \in \Mor(\calI)$
\end{itemize}
With this setup, the category $\overline{\scrR} \lRep$ coincides with the category $\scrR \lMod$ of left $\scrR$-modules as defined in \cite{SS2017}; see \cite[Theorem A.2]{DLLY}.
\end{exa}

\begin{exa}\label{tilder R}
Given a representation $\scrR$ of $\calI$ on {\sf Ring}, one  can define an $\calI^{\op}$-diagram $\widetilde{\scrR}$ of right module categories as follows:
\begin{itemize}
\item $\widetilde{\scrR}_i=\rMod R_i$, the category of right $R_i$-modules,
      for each $i \in \Ob(\calI\op)$;
\item $\widetilde{\scrR}_{\alpha\op}=-\otimes_{R_j}R_j: \rMod R_j \to \rMod R_i$ for each $\alpha\op: j \to i \in \Mor(\calI \op)$
\end{itemize}
Then $\widetilde{\scrR} \lRep$ coincides with the category $\rMod{\scrR}$ of right $\scrR$-modules; see \cite[Remark A.4]{DLLY}.
\end{exa}

Let $\scrD$ be a right exact $\calI$-diagram of Grothendieck categories such that each $\alpha \in \Mor(\calI)$, the functor $\scrD_\alpha$ preserves small coproducts. Then the colimit of a direct system $((M^x),(f^{yx}))$ of objects in $\scrD \lRep$ is defined componentwise; see \cite[1.12]{DLLY}. Dually, the limit of an inverse system $((M^x),(f^{xy}))$ of objects in $\scrD \lRep$ is also defined componentwise. The next result can be found in \cite[Theorem 2.8 and Proposition 2.10]{DLLY}.

\begin{thm} \label{Rep is grothendieck}
Let $\scrD$ be a right exact $\calI$-diagram of abelian categories such that each $\scrD_\alpha$ preserves small coproducts for $\alpha \in \Mor(\calI)$. If each $\scrD_i$ is a Grothendieck category (admitting enough projectives) for all $i \in \Ob(\calI)$, then so is $\scrD \lRep$. Moreover, if each $\scrD_i$ is locally finitely presented, then so is $\scrD \lRep$.
\end{thm}

\begin{setup*}
Throughout the paper, we use the notation $\scrD$ to denote a right exact $\calI$-diagram of Grothendieck categories, such that each $\scrD_\alpha$ preserves small coproducts for $\alpha \in \Mor(\calI)$.
\end{setup*}

In the following, we introduce some adjoint pairs that will be used frequently throughout the paper; see \cite[Section 3]{DLLY} for details.

\begin{ipg}\label{induction-func}
Let $\iota:\calJ \to \calI$ be a functor between skeletal small categories. It is easy to check that $\scrD\circ\iota$ is a right exact $\calJ$-diagram, where $(\scrD\circ\iota)_j = \scrD_{\iota(j)}$ for $j \in \Ob(\calJ)$ and $(\scrD \circ \iota)_\alpha = \scrD_{\iota(\alpha)}$ for $\alpha \in \Mor(\calJ)$. Then by \lemcite[2.2]{DLLY}, there exists an exact functor $\iota^*: \scrD\lRep \to (\scrD\circ\iota)\lRep,$ called the \textit{restriction functor} induced by $\iota$,  defined by $\iota^*(M)_j=M_{\iota(j)}$ and $\iota^*(M)_\alpha=M_{\iota(\alpha)}$ for each $M\in\scrD\lRep$. It follows from \prpcite[2.4]{DLLY} that $\iota^*$ admits a left adjoint functor $\iota_!: (\scrD \circ \iota) \lRep \to \scrD \lRep$, called the \textit{induction functor} induced by $\iota$.

Moreover, if  the $\calI$-diagram $\scrD$ admits enough right adjoints, then one can construct a right adjoint functor of $\iota^{\ast}$,  called the \emph{coinduction functor} induced by $\iota$, denoted by $\iota_{\ast}$; see \cite[Remark 2.5]{DLLY}.
\end{ipg}

\begin{ipg}\label{free-func}
Let $i\in \Ob(\calI)$. Define the functor $${\sf eva}^i: \scrD\lRep\to\scrD_i$$ which  sends a representation $M$ over $\scrD$ to its ``local" value $M_i$ in $\scrD_i$. By \corcite[2.6]{DLLY}, the  functor ${\sf eva}^i$ admits a left adjoint functor $\fre_i: \scrD_i\to\scrD\lRep$, which sends an object $M_i \in \scrD_i$ to $\fre_i(M_i)$ with
$(\fre_i(M_i))_j = \coprod_{\theta \in {\Hom[\calI]ij}} \scrD_{\theta} (M_i)$
for each $j \in \Ob(\calI)$.
Moreover, if the $\calI$-diagram $\scrD$ admits enough right adjoints, then one can construct a right adjoint functor of ${\sf eva}^i$, denoted by $\cofre_i$.
\end{ipg}

\section{Cofibrantly generated abelian model structures on $\scrD \lRep$}\label{sec3}
\noindent
In this section, we show that, under certain mild conditions,  a family of cofibrantly generated (hereditary) Hovey triples in $\scrD_i$'s induces a cofibrantly generated (hereditary) Hovey triple in $\scrD \lRep$.

\begin{nota}\label{nota od 4.2 and 4.3}
Throughout this section, let $(\calQ_i, \calW_i, \calR_i)$ be a Hovey triple in $\scrD_i$ for $i \in \Ob(\calI)$, and denote by
\begin{itemize}
\item $\widetilde{\calQ}_i$ (resp., $\widetilde{\calR}_i$) the subcategory $\calQ_i \cap \calW_i$ (resp., $\calW_i \cap \calR_i$) for $i \in \Ob(\calI)$;

\item $\calQ$ (resp., $\calW$, $\widetilde{\calQ}$, $\calR$, $\widetilde{\calR}$) the family $\{\calQ_i\}_{i \in \Ob(\calI)}$ (resp., $\{\calW_i\}_{i \in \Ob(\calI)}$, $\{\widetilde{\calQ}_i\}_{i \in \Ob(\calI)}$, $\{\calR_i\}_{i \in \Ob(\calI)}$, $\{\widetilde{\calR}_i\}_{i \in \Ob(\calI)}$) of subcategories;
\item $\scrD \lRep_{\calW} = \{M\in\scrD \lRep \,|\, M_i \in \calW_i \text{ for } i \in \Ob(\calI) \}$;

\item $\scrD \lRep_{\calR} = \{M\in\scrD \lRep \,|\,  M_i \in \calR_i \text{ for } i \in \Ob(\calI) \}$;

\item $\scrD \lRep_{\widetilde{\calR}} = \{M\in\scrD \lRep \,|\, M_i \in \widetilde{\calR}_i\text{ for } i \in \Ob(\calI) \}$.
\end{itemize}

By Theorem \ref{Hovey cor}, $(\Mon(\calQ_i), \overline{\sfW}_i, \Epi(\calR_i))$ forms an abelian model structure on $\scrD_i$ in which
\[
\overline{\sfW}_i = \{ w_i \,| \, w_i \text{ can be decomposed as } w_i = f_ic_i \text{ with }\ c_i \in \Mon(\widetilde{\calQ}_i) \text{ and } f_i \in \Epi(\widetilde{\calR}_i) \}.
\]

We always denote by
\begin{itemize}
\item $\overline{\sfW} = \{ \omega: M \to M' \in \scrD \lRep \,| \, \omega_i: M_i \to M'_i \in \overline{\sfW}_i \text{ for } i \in \Ob(\calI)\}$;

\item $\overline{\Epi}(\calR) = \{ \omega: M \to M' \in \scrD \lRep \,| \,      \omega_i: M_i \to M'_i \in \Epi(\calR_i) \text{ for } i \in \Ob(\calI) \}$.
\end{itemize}
\end{nota}

The next result can be proved using \lemcite[5.8]{Ho02}.

\begin{lem}\label{hovey}
For each $i \in \Ob(\calI)$, there are equalities
$$\Epi(\widetilde{\calR}_i) = \Epi(\calR_i) \cap \overline{\sfW}_i\ \ \mathrm{and}\ \ \Mon(\widetilde{\calQ}_i) = \Mon(\calQ_i) \cap \overline{\sfW}_i.$$
\end{lem}

The following two lemmas give equivalent characterizations for objects in the subcategories given in \ref{nota od 4.2 and 4.3}.

\begin{lem}\label{fibrant obs}
Let $M$ be an object in $\scrD \lRep$. Then the following statements hold:
\begin{prt}
\item $M \in \scrD \lRep_{\calW}$ if and only if $0 \to M \in \overline{\sfW}$ if and only if $M \to 0 \in \overline{\sfW}$.
\item $M \in \scrD \lRep_{\widetilde{\calR}}$ if and only if $M \to 0 \in \overline{\Epi}(\calR) \cap \overline{\sfW}$.
\item $M \in \scrD \lRep_{\calR}$ if and only if $M \to 0 \in \overline{\Epi}(\calR)$.
\end{prt}
\end{lem}

\begin{prf*}
(a). We only prove the first equivalence; the second one can be proved similarly. Suppose that $M \in \scrD \lRep_{\calW}$. We mention that $M_i \in \calW_i$ and $(\calQ_i, \calW_i, \calR_i)$ is a Hovey triple in $\scrD_i$. By Lemma \ref{model structor111}, there exists a short exact sequence $0 \to B_i \to A_i \overset{f_i}\longrightarrow M_i \to 0$ in $\scrD_i$ such that $A_i \in \widetilde{\calQ}_i$ and $B_i \in \widetilde{\calR}_i$. Since $f_i \in \Epi(\widetilde{\calR}_i)$ and $0 \to A_i \in \Mon(\widetilde{\calQ}_i)$, it follows that $0 \to M_i \in \overline{\sfW}_i$, as desired.

Conversely, if $0 \to M \in \overline{\sfW}$, then $0 \to M_i \in \overline{\sfW}_i$ and hence can be factored as $0 \to A'_i \overset{h_i} \longrightarrow M_i$ with $0 \to A'_i \in \Mon(\widetilde{\calQ}_i)$ and $h_i \in \Epi(\widetilde{\calR}_i)$. Consequently, $A'_i \in \widetilde{\calQ}_i = \calQ_i \cap \calW_i$ and there exists a short sequence $0 \to K_i \to A'_i \overset{h_i} \longrightarrow M_i \to 0$ in $\scrD_i$ with $K_i \in \widetilde{\calR}_i = \calR_i \cap \calW_i$. But $\calW_i$ is thick, so $M_i \in \calW_i$.

(b). By Lemma \ref{hovey}, one has $\Epi(\widetilde{\calR}_i) = \Epi(\calR_i) \cap \overline{\sfW}_i$. Thus one gets that $M \in \scrD \lRep_{\widetilde{\calR}}$ if and only if $M_i \to 0 \in \Epi(\widetilde{\calR}_i)$ for each $i \in \Ob(\calI)$, if and only if $M_i \to 0 \in \Epi(\calR_i) \cap \overline{\sfW}_i$ for $i \in \Ob(\calI)$, that is, $M \to 0 \in \overline{\Epi}(\calR) \cap \overline{\sfW}$.

(c). One gets that $M \in \scrD \lRep_{\calR}$ if and only if $M_i \to 0 \in \Epi(\calR_i)$ for each $i \in \Ob(\calI)$, that is, $M \to 0 \in \overline{\Epi}(\calR)$.
\end{prf*}

\begin{lem}\label{abelian epi}
Let $\omega$ be a morphism in $\scrD \lRep$. Then the following statements hold:
\begin{prt}
\item $\omega \in \overline{\Epi}(\calR) \cap \overline{\sfW}$ if and only if it is an epimorphism and $\ker(\omega) \to 0 \in \overline{\Epi}(\calR) \cap \overline{\sfW}$.

\item $\omega \in \overline{\Epi}(\calR)$ if and only if it is an epimorphism and $\ker(\omega) \to 0 \in \overline{\Epi}(\calR)$.
\end{prt}
\end{lem}

\begin{prf*}
We only give a proof for statement (a); one can prove (b) similarly. By Lemma \ref{hovey}, one has $\Epi(\calR_i) \cap \overline{\sfW}_i = \Epi(\widetilde{\calR}_i)$ for each $i \in \Ob(\calI)$. Therefore, $\omega \in \overline{\Epi}(\calR) \cap \overline{\sfW}$ if and only if $\omega_i \in \Epi(\widetilde{\calR}_i)$ for each $i \in \Ob(\calI)$. The conclusion then follows from Lemma \ref{fibrant obs}(c).
\end{prf*}

\begin{dfn} \label{compatibility}
Let $\mathcal{S} = \{\mathcal{S}_i\}_{i \in \Ob(\calI)}$ be a family with each $\mathcal{S}_i$ a subcategory of $\scrD_i$. We say that $\mathcal{S}$ is \emph{compatible with  respect to} $\scrD$ if $\scrD_\alpha (\mathcal{S}_{i}) \subseteq \mathcal{S}_{j}$ for each $\alpha: i\to j$ in $\Mor(\calI)$.
\end{dfn}

\begin{rmk}
Actually, by comparing the above definition to the one of subdiagrams (see \cite[Definition 1.4]{DLLY}), one can easily see that $\mathcal{S} = \{\mathcal{S}_i\}_{i \in \Ob(\calI)}$ is compatible with  respect to $\scrD$ if and only if it is a subdiagram of $\scrD$. In this case, we denote this subdiagram by $\mathcal{S}$ by abuse of notation. In particular, if an object $M$ in $\scrD \lRep$ satisfies the condition that $M_i \in \mathcal{S}_i$ for each $i \in \Ob(\calI)$, then $M$ is actually a representation over the subdiagram $\mathcal{S}$. Consequently, one has the following obvious identification
\[
\scrD \lRep_{\mathcal{S}} = \{M \in \scrD \lRep \, \mid \, M_i \in \mathcal{S}_i \text{ for } i \in \Ob(\calI) \} = \mathcal{S} \lRep.
\]
\end{rmk}

By Lemma \ref{hovey}, there is an equality $\Mon(\widetilde{\calQ}_i) = \Mon(\calQ_i) \cap \overline{\sfW}_i$ for each $i \in \Ob(\calI)$. This identity allows us to deduce the following result immediately.

\begin{lem}\label{keep c and tc}
Suppose that $\scrD$ is exact. If both $\calQ$ and $\widetilde{\calQ}$ are compatible with respect to $\scrD$, then $\scrD_\alpha(\Mon(\calQ_i)) \subseteq \Mon(\calQ_j)$ and $\scrD_\alpha(\Mon(\calQ_i) \cap \overline{\sfW}_i) \subseteq \Mon(\calQ_j) \cap \overline{\sfW}_j$ for all $\alpha: i \to j$ in $\Mor(\calI)$.
\end{lem}

Now assume that the Hovey triple $(\calQ_i, \calW_i, \calR_i)$ is cofibrantly generated for each $i \in \Ob(\calI)$, which means that the associated abelian model structure $(\Mon(\calQ_i), \overline{\sfW}_i, \Epi(\calR_i))$ on $\scrD_i$ is cofibrantly generated. Denote by
\begin{itemize}
\item $I_i$ the set of generating cofibrations of $(\Mon(\calQ_i), \, \overline{\sfW}_i, \, \Epi(\calR_i))$,

\item $J_i$ the set of generating trivial cofibrations of $(\Mon(\calQ_i), \, \overline{\sfW}_i, \, \Epi(\calR_i))$,

\item $\fre_\bullet(I_\bullet) = \{ \fre_i(f_i) \,|\,  f_i \in I_i \text{ and } i \in \Ob(\calI) \}$,

\item $\fre_\bullet(J_\bullet) = \{ \fre_i(g_i) \,|\, g_i \in J_i \text{ and } i \in \Ob(\calI) \}$.
\end{itemize}

\begin{prp}\label{HR g MODEL}
Suppose that $\scrD$ is exact.
If $(\calQ_i, \calW_i, \calR_i)$ is a cofibrantly generated Hovey triple in $\scrD_i$  for $i \in \Ob(\calI)$, and both $\calQ$ and $\widetilde{\calQ}$ are compatible with respect to $\scrD$, then
$$({^\Box(\fre_\bullet(I_\bullet)^\Box)}, \, \overline{\sfW}, \, \overline{\Epi}(\calR))$$
is a cofibrantly generated abelian model structure on $\scrD \lRep$.
\end{prp}

\begin{prf*}
We mention that $(\Mon(\calQ_i), \sfW_i, \Epi(\calR_i))$ is a cofibrantly generated abelian model structure on $\scrD_i$ for each $i \in \Ob(\calI)$ with $I_i$ (resp., $J_i$) the set of generating cofibrations (resp., trivial cofibrations). Then $\scrD_\alpha$ preserves cofibrations and trivial cofibrations for $\alpha \in \Mor(\calI)$ by Lemma \ref{keep c and tc}. Note that  both $\fre_\bullet(I_\bullet)$ and $\fre_\bullet(J_\bullet)$ are indeed sets, and the triplet $$({^\Box(\overline{\sfW}\cap\fre_\bullet(J_\bullet)^\Box)}, \, \overline{\sfW}, \, \fre_\bullet(J_\bullet)^\Box)$$
defines a $g$-structure in the sense of H\"{u}ttemann and R\"{ondigs}; see \cite[Definition 3.4.2]{HR08}\footnote{Indeed, in \cite{HR08}, $\mathcal{B}$ is an $\mathcal{I}$-bundle of model categories, which is actually a diagram of model categories admitting enough right adjoints satisfying the condition that $\mathcal{B}_{\alpha}$ preserves cofibrations and trivial cofibrations for $\alpha \in \Mor(\calI)$, and the category ${\bf Tw}(\calI, \mathcal{B})$ of twisted diagrams is actually the category $\mathcal{B}\lRep$ in our sense. We mention that the results in \cite{HR08} that we used in Sections \ref{sec3} and \ref{abelian model structures and cotorsion pairs} still hold without the assumption that $\mathcal{B}$ admits enough right adjoints.}. Thus it follows from \thmcite[3.4.5]{HR08} that the above triplet is a model structure on $\scrD \lRep$ cofibrantly generated by $\fre_\bullet(I_\bullet)$ and $\fre_\bullet(J_\bullet)$, and so one has $\fre_\bullet(I_\bullet)^\Box=\overline{\sfW}\cap\fre_\bullet(J_\bullet)^\Box$. Note that by \cite[Lemma 3.4.3]{HR08} that a morphism $f: M \to N$ in $\scrD\lRep$ has the right lifting property with respect to $\fre_\bullet(J_\bullet)$ if and only if for each object $i$, $f_i$ is a fibration in $\scrD_i$, that is, $f_i\in\Epi(\calR_i)$. This yields that $\fre_\bullet(J_\bullet)^\Box = \overline{\Epi}(\calR)$. Hence the model structure can be rewritten as
$$({^\Box(\fre_\bullet(I_\bullet)^\Box)}, \, \overline{\sfW}, \, \overline{\Epi}(\calR)).$$
It remains to show that it is abelian. By Lemma \ref{abelian epi}, it suffices to show that any $\omega \in {^\Box(\fre_\bullet(I_\bullet)^\Box)}$ is a monomorphism. Since $I_i \subseteq \Mon(\calQ_i)$, each morphism in $I_i$ is a monomorphism. Therefore, all morphisms in $\fre_\bullet(I_\bullet)$ are monomorphisms since $\scrD$ is exact. By the general theory of cofibrantly generated model structures, we see that $\omega$ is a retract of a transfinite composition of pushouts of monomorphisms. Thus $\omega$ is a monomorphism as desired; see \prpcite[A.6.(2)]{SS2011}.
\end{prf*}

We are now ready to present the main result of this section.

\begin{thm}\label{ht induce ht 2}
Suppose that $\scrD$ is exact. If $(\calQ_i, \calW_i, \calR_i)$ is a cofibrantly generated (hereditary) Hovey triple in $\scrD_i$  for $i \in \Ob(\calI)$, and both $\calQ$ and $\widetilde{\calQ}$ are compatible with respect to $\scrD$, then
\[
(^\perp\scrD \lRep_{\widetilde{\calR}}, \, \scrD \lRep_{\calW}, \, \scrD \lRep_{\calR})
\]
is a cofibrantly generated (hereditary) Hovey triple in $\scrD \lRep$.
\end{thm}

\begin{prf*}
Since
$({^\Box(\fre_\bullet(I_\bullet)^\Box)}, \, \overline{\sfW}, \, \overline{\Epi}(\calR))$ forms a cofibrantly generated abelian model structure on $\scrD \lRep$ by Proposition \ref{HR g MODEL}, it follows that the associated subcategories of cofibrant, trivial and fibrant objects form a cofibrantly generated Hovey triple in $\scrD \lRep$. By Lemma \ref{fibrant obs}(a), the subcategory of trivial objects is $\scrD \lRep_{\calW}$. It follows from Lemma \ref{fibrant obs}(c) that the subcategory of fibrant objects is $\scrD \lRep_{\calR}$. Thus the subcategory of cofibrant objects is
$^\perp(\scrD \lRep_{\calW}\cap\scrD \lRep_{\calR}) = {^\perp(\scrD \lRep_{\widetilde{\calR}})}$.
Consequently, the triplet
\[
({^\perp\scrD \lRep_{\widetilde{\calR}}}, \, \scrD \lRep_{\calW}, \, \scrD \lRep_{\calR})
\]
forms a cofibrantly generated Hovey triple in $\scrD \lRep$.

Suppose in addition that each Hovey triple $(\calQ_i, \calW_i, \calR_i)$ is hereditary. Then both $\widetilde{\calR}_i$ and $\calR_i$ are closed under  cokernels of monomorphisms for $i \in \Ob(\calI)$. To show the hereditary property of the above Hovey triple, it suffices to prove the hereditary property of the complete cotorsion pairs
\begin{center}
$({^\perp\scrD \lRep_{\widetilde{\calR}}} \cap \scrD \lRep_{\calW}, \, \scrD \lRep_{\calR})$
and
$({^\perp\scrD \lRep_{\widetilde{\calR}}}, \, \scrD \lRep_{\widetilde{\calR}})$.
\end{center}
It is clear that both $\scrD \lRep_{\widetilde{\calR}}$ and $\scrD \lRep_{\calR}$ are closed under taking cokernels of monomorphisms. Therefore, the above two cotorsion pairs are coresolving, and hence hereditary by Lemma \ref{test hereditary}.
\end{prf*}

\section{Induced abelian model structures on $\scrD \lRep$}
\label{abelian model structures and cotorsion pairs}
\noindent
The abelian model structure on $\scrD \lRep$ constructed in the previous section  works for any skeletal small index category $\calI$, but it has two limitations: it requires the given family of abelian model structures indexed by objects in $\calI$ to be cofibrantly generated, and it does not provide an explicit description of the cofibrant objects. In this section, we focus a special kind of index categories $\calI$, namely left rooted quivers (viewed as categories in a natural way). It turns out that for such categories, we are able to provide  an explicit description of cofibrant objects, which are closely related to monomorphism categories studied in \cite{Z11}.

\begin{bfhpg}[Rooted categories]\label{construct}
Suppose that $\calI$ is a partially ordered category, that is, where the relation $\preccurlyeq$ on $\Mor(\calI)$, defined by setting $i \preccurlyeq j$ if $\Hom[\calI]i{j} \neq \emptyset$ is a partial order. Define a transfinite sequence $\{V_{\chi}\}_{\chi \, \mathrm{ordinal}}$ of subsets of $\Ob(\calI)$ as follows:
\begin{itemize}
\item for the first ordinal $\chi = 0$, set $V_0 = \emptyset$;

\item for a successor ordinal $\chi+1$, set
 \[
V_{\chi+1} = \left \{ i \in \Ob(\calI) \:
      \left|
        \begin{array}{c}
          i \text{ is not the target of any } \alpha \in \Mor(\calI) \\
          \text{with source } s(\alpha)\neq i \text{ and } s(\alpha) \notin \cup_{\mu \leqslant \chi} V_{\mu}
        \end{array}
      \right.
    \right\};
\]

\item for a limit ordinal $\chi$, set $V_{\chi} = \cup_{\mu < \chi} V_{\mu}$.
\end{itemize}

Following \cite[Definition 3.4]{DLLY}, we say  that a partially ordered category $\calI$ is \emph{left rooted} if there exists an ordinal $\zeta$ such that $V_{\zeta} = \Ob(\calI)$.
We say that $\calI$ is \emph{right rooted} if $\calI\op$ is left rooted.
\end{bfhpg}

\begin{exa}\label{rooted}
Let $\calI$ be a quiver (viewed as a category in a natural way).
There exists a transfinite sequence $\{V_{\alpha}\}$ of subsets of $\Ob(\calI)$ as follows:
\begin{prt}
\item[$\bullet$] For the first ordinal $\alpha=0$ set $V_0=\emptyset$, for a successor ordinal $\alpha+1$ set
$$V_{\alpha+1}=\{i\in \Ob(\calI)~|~i\ \text{is not the target of any arrow}\ a\in \calI\ \text{with}\ s(a)\notin\cup_{\beta\leq\alpha}V_{\beta}\},$$
and for a limit ordinal $\alpha$ set $V_{\alpha}=\cup_{\beta<\alpha}V_{\beta}$.
\end{prt}
Clearly, the sets form a chain
$V_1\subseteq V_2\subseteq\cdots\subseteq \Ob(\calI)$. Recall from \dfncite[3.5]{EOT04} that a quiver $\calI$ is called \emph{left rooted} if there exists an ordinal $\lambda$ such that $V_{\lambda}=\Ob(\calI)$. By \prpcite[3.6]{EOT04}, a quiver $\calI$ is left rooted if and only if it has no infinite sequence of arrows of the form $\cdots \to \bullet\to \bullet \to \bullet$ (not necessarily different), and so there is no loop or oriented cycle in a left rooted quiver. We mention that a quiver $\calI$ without loops or oriented cycles is left rooted if and only if, when viewed as a category, it is a  left rooted category; see \cite[Remark 3.5]{DLLY}.
\end{exa}

The following definition is taken from \cite[Definition 5.1.1]{Ho99}.

\begin{dfn} \label{dfn of direct categories}
A skeletal small category $\calI$ is called a \textit{direct category} if there exists a functor $F: \calI \to \zeta$, where $\zeta$ is an ordinal (viewed as a category in a natural way) such that $F$ sends non-identity morphisms in $\calI$ to non-identity morphisms in $\zeta$. We say that $\calI$ is an \emph{inverse} category if $\calI\op$ is a direct category.
\end{dfn}

\begin{rmk}\label{direct-rooted}
It follows from \cite[Proposition 3.7]{DLLY} that $\calI$ is direct if and only if it is left rooted and locally trivial (that is, the set $\AU i$ of endomorphisms on $i$ contains only the identity morphism for all $i \in \Ob(\calI)$).
\end{rmk}

\begin{ipg}\label{latching functor}
Let $\calI$ be a direct category.
Fix $i \in \Ob(\calI)$ and denote by $\calI_{\prec i}$ the full subcategory consisting of objects $j \in \Ob(\calI)$ such that $j \prec i$. For each $M \in \scrD \lRep$, consider the component of the counit $\iota_!\iota^*(M) \to M$ of the adjoint pair $(\iota_!, \iota^*)$ at $M$, where $\iota^*$ is the restriction functor induced by the natural embedding functor $\iota: \calI_{\prec i} \to \calI$; see \ref{induction-func}. Restricting to the object $i$, we obtain an object $\sfL_i(M) = (\iota_! \iota^*(M))_i$ as well as a natural morphism $\sfL_i(M) \to M_i$ in $\scrD_i$. Indeed, one has
$$\sfL_i(M)=\colim_{\substack{\theta: \, h \to i \in \Mor(\calI) \\ h \neq i}}\scrD_\alpha(M_{h})=\colim_{\alpha \in \calP_i(\bullet, i)}\scrD_\alpha(M_{s(\alpha)}),$$
where $\calP_i = \Mor(\calI) \backslash \AU i$ is a prime ideal of $\calI$ in the sense of \cite[Subsection 2.3]{DLLY} as $\calI$ is a partially ordered category. Here $\AU i$ is the set of endomorphisms on $i$. It is easy to see that $\sfL_i$ is precisely the \emph{latching functor} described in \cite{HR08}, and the natural morphism $\sfL_i(M) \to M_i$ is actually the morphism $\varphi_i^M: \colim_{\alpha \in \calP_i(\bullet, i)} \scrD_\alpha(M_{s(\alpha)})\to M_i$ given in \cite[(2.12.1)]{DLLY}.

\begin{dfn} \label{s and phi}
Given a direct category $\calI$ and a family $\calX = \{ \calX_i \}_{i \in \Ob(\calI)}$ with each $\calX_i$ a full subcategory of $\scrD_i$, define a subcategory of $\scrD \lRep$:
\[
\Phi(\calX)= \{ M \in \scrD \lRep \mid \varphi_i^M \text{ is a monomorphism and } \coker(\varphi_i^M) \in {\calX}_i \text{ for each } i \in \Ob(\calI) \}.
\]
In particular, one has
\[
\Phi(\scrD) = \{ X \in \scrD \lRep \mid \varphi_i^X \text{ is a monomorphism for each } i \in \Ob(\calI) \}.
\]
\end{dfn}

\begin{rmk}\label{exa4.23}
If $\calI$ is a left rooted quiver (viewed as a
category in a natural way), then it is direct (see Remark \ref{direct-rooted}), and the colimit appearing in the definition of $\sfL_i$ is actually a coproduct, that is, $$\sfL_i(M)=\coprod_{\alpha \in \calI(\bullet, i)}\scrD_\alpha(M_{s(\alpha)}),$$
where $\calI(\bullet, i)$ denotes the set of all arrows in $\calI$ with $i$ the target. In this case, the morphism $\varphi_i^M$ is from $\coprod_{\alpha \in \calI(\bullet, i)}\scrD_\alpha(M_{s(\alpha)})$ to $M_i$.
In some subsequent results in this section, we often assume that $\calI$ is a left rooted quiver; the main obstacle forcing us to work with left rooted quivers rather than arbitrary direct categories is that the colimit functor (not necessarily a filtered colimit) is not exact in general.
\end{rmk}

For any morphism $\omega : M \to N$ in $\scrD \lRep$, by considering the following commutative diagram where the inner square is a pushout, we have a natural morphism $\rho_i$:
\begin{equation} \label{latch pushout}
\xymatrix{
  \sfL_i(M)  \ar[d]_{\sfL_i(\omega)} \ar[r]^{\varphi_i^M} & M_i \ar[d]_{\theta_i} \ar@/^/[ddr]^{\omega_i}  \\
  \sfL_i(N) \ar[r]^-{\delta_i} \ar@/_/[drr]_{\varphi_i^{N}}  & M_i \sqcup_{\sfL_i(M)} \sfL_i(N) \ar@{.>}[dr]|-{\rho_i} \\
                &               & N_i.}
\end{equation}
Set $\overline{\POMon}(\calQ) = \{ \omega: M \to N \in \scrD \lRep \,|\,  \rho_i: M_i \sqcup_{\sfL_i (M)} \sfL_i (N) \to N_i \in \Mon(\calQ_i) \text{ for all } i \in \Ob(\calI)\}$.
\end{ipg}

\begin{lem}\label{confibration mono1}
Suppose that $\calI$ is a left rooted quiver, and $\scrD$ is exact.
Then any morphism in $\overline{\POMon}(\calQ)$ is a monomorphism.
\end{lem}

\begin{prf*}
Let $\{V_{\chi}\}_{\chi\,{\rm ordinal}}$ be the transfinite sequence of subsets of $\Ob(\calI)$ defined in \ref{construct}. Since $\calI$ is a left rooted quiver, there exists an ordinal $\zeta$ such that $\Ob(\calI) = V_{\zeta}$. Take $\omega: M \to N \in \overline{\POMon}(\calQ)$. We will use the transfinite induction to show that $\omega_i: M_i \to N_i$ is a monomorphism for all ordinals $\chi \leqslant \zeta$ and all $i \in V_{\chi}$. Consider the commutative diagram (\ref{latch pushout}). Then $\rho_i$ is a monomorphism for $i \in \Ob(\calI)$ as $\omega\in\overline{\POMon}(\calQ)$.

If $i \in V_1$, the set of minimal objects with respect to the partial order $\preccurlyeq$ defining the partially ordered structure of $\calI$, then $\sfL_i (M) = 0 = \sfL_i (N)$ as $\calI(\bullet, i) = \emptyset$, so $\omega_i = \rho_i$ as $\theta_i = \id_{M_i}$, and hence $\omega_i$ is a monomorphism.

For $\chi > 1$, we have two cases:

(1) If $\chi$ is a successor ordinal and $i \in V_{\chi}$, then $j \in V_{\chi-1}$ for all $j \prec i$ (see \cite[Remark 3.3]{DLLY}). By the induction hypothesis, all $\omega_j$ are monomorphisms. Since $\scrD$ is exact, $\scrD_{\alpha} (\omega_j): \scrD_{\alpha} (M_j) \to \scrD_{\alpha} (N_j)$ is also a monomorphism for $j \prec i$. It follows that $\sfL_i(\omega)$ is a monomorphism as well; see Remark \ref{exa4.23}. Consequently, $\theta_i$ is a monomorphism since the inner square is a pushout. Thus $\omega_i = \rho_i\theta_i$ is a monomorphism.

(2) If $\chi \leqslant \zeta$ is a limit ordinal, then the conclusion is clearly true for $\chi$ because in this case $V_{\chi} = \cup_{\mu < \chi} V_{\mu}$ and the conclusion holds for all ordinals $\mu < \chi$.
\end{prf*}

\begin{prp}\label{Rondigs}
Suppose that $\calI$ is a left rooted quiver, and $\scrD$ is exact. If $(\calQ_i, \calW_i, \calR_i)$ is a Hovey triple in $\scrD_i$ for $i \in \Ob(\calI)$, and both $\calQ$ and $\widetilde{\calQ}$ are compatible with respect to $\scrD$, then
\[
(\overline{\POMon}(\calQ), \, \overline{\sfW}, \, \overline{\Epi}(\calR))
\]
is an abelian model structure on $\scrD \lRep$.
\end{prp}

\begin{prf*}
We mention that $(\Mon(\calQ_i), \sfW_i, \Epi(\calR_i))$ is an abelian model structure on $\scrD_i$ for each $i \in \Ob(\calI)$ by Theorem \ref{Hovey cor}. Since $\calQ$ and $\widetilde{\calQ}$ are compatible with respect to $\scrD$, Lemma \ref{keep c and tc} tells us that $\scrD_\alpha$ preserve cofibrations and acyclic cofibrations for $\alpha \in \Mor(\calI)$. It is easy to see that the triplet
$$(\overline{\POMon}(\calQ), \, \overline{\sfW}, \, \overline{\Epi}(\calR))$$
is the $c$-structure in the sense of H\"{u}ttemann and R\"{ondigs}; see \cite[Definition 3.2.6]{HR08}. Thus it follows from \thmcite[3.2.13(1)]{HR08} that the above triplet is a model structure on $\scrD \lRep$, and is furthermore abelian by Lemmas \ref{abelian epi} and \ref{confibration mono1} as well as Definition \ref{df of abelian ms}.
\end{prf*}

\begin{rmk}
\thmcite[3.2.13]{HR08} tells us that a compatible family of model structures can be amalgamated to a model structure of the above form on $\scrD \lRep$. Thus one may wonder if the family of model structures are abelian model structures, then their amalgamation is also abelian. This may be not true in general (though we do not have a counterexample at hand). The above proposition gives a sufficient criterion such that the amalgamation is indeed abelian.
\end{rmk}

The following result gives a description of cofibrant objects in the abelian model structure described in Proposition \ref{Rondigs}.

\begin{lem}\label{cofibrant obs phi}
Suppose that $\calI$ is a left rooted quiver. Then an object $N$ in $\scrD \lRep$ is contained in $\Phi(\calQ)$ if and only if the morphism $0 \to N$ lies in $\overline{\POMon}(\calQ)$.
\end{lem}

\begin{prf*}
Recall that $N \in \Phi(\calQ)$ if and only if $\varphi_i^{N}: \sfL_i(N) \to N_i$ is a monomorphism and $\coker(\varphi_i^{N})\in\calQ_i$ for each $i \in \Ob(\calI)$. Take $M = 0$ in the commutative diagram (\ref{latch pushout}). Then $\delta_i = \id_{\sfL_i(N)}$, $\varphi_i^{N} = \rho_i$, so $N \in \Phi(\calQ)$ if and only if $\rho_i$ is a monomorphism and $\coker(\rho_i) \in \calQ_i$ for each $i \in \Ob(\calI)$, which is equivalent to saying that $0 \to N$ is contained in $\overline{\POMon}(\calQ)$.
\end{prf*}

We are now ready to give the main result of this section.

\begin{thm}\label{ht induce ht 1}
Suppose that $\calI$ is a left rooted quiver, and $\scrD$ is exact. If $(\calQ_i, \calW_i, \calR_i)$ is a Hovey triple in $\scrD_i$ for $i \in \Ob(\calI)$, and both $\calQ$ and $\widetilde{\calQ}$ are compatible with respect to $\scrD$. Then
\[
(\Phi(\calQ), \, \scrD \lRep_{\calW}, \, \scrD \lRep_{\calR})
\]
forms a Hovey triple in $\scrD \lRep$.
If furthermore, each Hovey triple $(\calQ_i, \calW_i, \calR_i)$ is hereditary for $i \in \Ob(\calI)$, then the above Hovey triple in $\scrD \lRep$ is hereditary as well.
\end{thm}

\begin{prf*}
Since $(\overline{\POMon}(\calQ), \, \overline{\sfW}, \, \overline{\Epi}(\calR))$ is an abelian model structure on $\scrD \lRep$ by Proposition \ref{Rondigs}, we only need to recognize subcategories of cofibrant, trivial and fibrant objects. The subcategory of cofibrant objects is $\Phi(\calQ)$ by Lemma \ref{cofibrant obs phi}, the subcategory of trivial objects is $\scrD \lRep_{\calW}$ by Lemma \ref{fibrant obs}(a), and the subcategory of fibrant objects is $\scrD \lRep_{\calR}$ by Lemma \ref{fibrant obs}(c). The first statement then follows. The second statement can be established using a similar argument as in the proof of Theorem \ref{ht induce ht 2}.
\end{prf*}

\begin{exa}
Let $\calA$ be a Grothendieck category, and let $\scrD$ be a trivial $\calI$-diagram of $\calA$. Then one has $\scrD \lRep={\sf Rep}(\calI, \calA)$, the category of representations of $\calI$ with values in $\calA$. In this case, Theorem \ref{ht induce ht 1} can be rewritten as: If $\calI$ is a left rooted quiver and $\calA$ is a Grothendieck category, then any Hovey triple $(\calQ, \calW, \calR)$ in $\calA$ induces a Hovey triple
$(\Phi(\calQ), \, {\sf Rep}(\calI, \calW), \, {\sf Rep}(\calI, \calR))$ in ${\sf Rep}(\calI, \calA)$. This fact improves \cite[Theorem B]{DELO} by removing the unnecessary condition that the Hovey triple
$(\calQ, \calW, \calR)$ is hereditary, which is essential in their proof since Hovey's correspondence was used.
\end{exa}

A careful reader may observe the following subtle fact. We use the subcategory $\calQ_i$ of cofibrant objects in $\scrD_i$ to construct the subcategory $\Phi(\calQ)$ of cofibrant objects in $\scrD \lRep$. We can also construct a subcategory $\Phi(\widetilde{\calQ})$ in $\scrD \lRep$ using the subcategory $\calQ_i \cap \calW_i$ of trivial confibrant objects in $\scrD_i$. A natural question is: under what conditions is $\Phi(\widetilde{\calQ})$ exactly the subcategory of trivial cofibrant objects in $\scrD \lRep$, that is, $\Phi(\widetilde{\calQ}) = \Phi(\calQ) \cap \scrD \lRep_{\calW}$\,? An answer is given in the rest of this section.

\begin{lem} \label{X(i) in X}
Let $\mathcal{S} = \{\mathcal{S}_i\}_{i \in \Ob(\calI)}$ be a family with each $\mathcal{S}_i$ a subcategory of $\scrD_i$, and suppose that $\mathcal{S}$ is compatible with respect to $\scrD$. If $\calI$ is a left rooted quiver, and $\mathcal{S}_i$ is closed under extensions and small coproducts for each $i \in \Ob(\calI)$, then $\Phi(\mathcal{S}) \subseteq \scrD \lRep_{\mathcal{S}}$.
\end{lem}

\begin{prf*}
Let $\{V_{\chi}\}_{\chi\,{\rm ordinal}}$ be the transfinite sequence of subsets of $\Ob(\calI)$ defined in \ref{construct}. Then one has $\Ob(\calI) = V_{\zeta}$ for a certain ordinal $\zeta$. Take $S \in {\Phi}({\mathcal{S}})$, we want to show that $S_i \in \mathcal{S}_i$ for all ordinals ${\chi}$ and all objects $i \in V_{\chi}$. This is trivially true for $V_0 = \emptyset$.

For $\chi \geqslant 1$, we have two cases:

(1) If $\chi$ is a successor ordinal, then for $i \in V_{\chi}$, by the definition of $\Phi(\mathcal{S})$, there is a short exact sequence
\[
0 \to \coprod_{\alpha \in \calI(\bullet, i)} \scrD_\alpha(S_{s(\alpha)}) \overset{\varphi^S_i} {\longrightarrow} S_i \to \coker(\varphi^S_i) \to 0
\]
in $\scrD_i$ with $\coker(\varphi^S_i)\in {\mathcal{S}_i}$. Since $\mathcal{S}_i$ is closed under coproducts and extensions, it suffices to show that each $\scrD_\alpha(S_{s(\alpha)})$ is contained in $\mathcal{S}_i$. But this is obvious. Indeed, by \cite[Remark 3.3]{DLLY}, $s(\alpha) \in V_{\chi-1}$, so $S_{s(\alpha)} \in \mathcal{S}_{s(\alpha)}$ by the induction hypothesis, and hence $\scrD_\alpha(S_{s(\alpha)}) \in \mathcal{S}_i$ since $\mathcal{S}$ is compatible with respect to $\scrD$.

(2) If $\chi \leqslant \zeta$ is a limit ordinal, then the conclusion clearly holds for $\chi$ because in this case $V_{\chi} = \cup_{\mu < \chi} V_{\mu}$ and the conclusion holds for all ordinals $\mu < \chi$.

The conclusion then follows by the transfinite induction.
\end{prf*}

\begin{prp}\label{com and her cp for q tilte}
Suppose that $\calI$ is a left rooted quiver, and $\scrD$ is exact.
If $(\calQ_i, \calW_i, \calR_i)$ is a Hovey triple in $\scrD_i$ for $i \in \Ob(\calI)$, and both $\calQ$ and $\widetilde{\calQ}$ are compatible with respect to $\scrD$. Then $\Phi(\widetilde{\calQ}) = \Phi(\calQ) \cap \scrD \lRep_{\calW}$.
\end{prp}

\begin{prf*}
The inclusion $\Phi(\widetilde{\calQ}) \subseteq \Phi(\calQ)$ hold trivially. Furthermore, by the previous lemma, one has
\[
\Phi(\widetilde{\calQ}) \subseteq \scrD \lRep_{\widetilde{\calQ}} \subseteq \scrD \lRep_{\calW}.
\]
From these two inclusion we deduce that $\Phi(\widetilde{\calQ}) \subseteq \Phi(\calQ) \cap \scrD \lRep_{\calW}$.

Conversely, taking an arbitrary $M \in \Phi(\calQ)\cap \scrD \lRep_{\calW}$, we want to show $M \in \Phi(\widetilde{\calQ})$. For each $i \in \Ob(\calI)$, there is a short exact sequence
\[
0 \to \coprod_{\alpha \in \calI(\bullet, i)} \scrD_\alpha(M_{s(\alpha)}) \overset{\varphi^M_i} {\longrightarrow} M_i \to \coker(\varphi^M_i) \to 0
\]
in $\scrD_i$ with $\coker(\varphi^M_i)\in \calQ_i$. By the definition of $\Phi(\widetilde{\calQ})$, it suffices to check that $\coker(\varphi^M_i) \in \calW_i$.
Indeed, since $M \in \Phi(\calQ)$, it follows from Lemma \ref{X(i) in X} that $M_i \in \calQ_i$ for all $i \in \Ob(\calI)$. Note that $M_i \in \calW_i$ as well by assumption, so $M_i \in \widetilde{\calQ}_i$ for $i \in \Ob(\calI)$. But $\widetilde{\calQ}$ is compatible with respect to $\scrD$ and $\widetilde{\calQ}_i$ is closed under small coproducts. Thus one has
$\coprod_{\alpha \in \calI(\bullet, i)} \scrD_\alpha(M_{s(\alpha)}) \in \widetilde{\calQ}_i \subseteq \calW_i$.
By the 2-out-of-3 property, $\coker(\varphi^M_i)\in \calW_i$ as desired.
\end{prf*}

\section{Gorenstein injective model structure on $\widetilde{\scrR} \lRep$}
\label{representations of small categories}%%%%%%%%%%
\noindent
In the rest of the paper, we turn our attention to two specific categories $\widetilde{\scrR} \lRep$ and $\overline{\scrR} \lRep$ (see Examples \ref{tilder R} and \ref{bar R}), and present some applications of results developed in Section \ref{abelian model structures and cotorsion pairs}. We first construct the Gorenstein injective model structure on $\widetilde{\scrR} \lRep$ in this section, where $\widetilde{\scrR}$ is the $\calI\op$-diagram of right module categories induced by a representation $\scrR$ of $\calI$ on {\sf Ring} with $\widetilde{\scrR}_i=\rMod R_i$ for $i \in \Ob(\calI\op)$ and $\widetilde{\scrR}_{\alpha\op}=-\otimes_{R_j}R_j: \rMod R_j \to \rMod R_i$ for $\alpha\op: j \to i \in \Mor(\calI \op)$ (see Example \ref{tilder R}). We then provide a characterization of Gorenstein injective objects in $\widetilde{\scrR} \lRep$. These results are not only of independent interest, but also play a key role in constructing the Gorenstein flat model structure and and in characterizing Gorenstein flat objects in the category $\overline{\scrR} \lRep$.

\begin{ipg}
For an arbitrary associative ring $A$, recall from Enochs and Jenda \cite{EEnOJn95b} that a right $A$-module $N$ is called \emph{Gorenstein injective} if there is an exact sequence
\[
\mathbb{I}: \, \cdots \to I^{-1} \to I^0 \to I^1 \to \cdots
\]
of injective right $A$-modules such that $N \cong \ker{(I^0 \to I^1)}$ and the sequence $\mathbb{I}$ remains exact after applying the functor $\Hom[A]{E}{-}$ for every injective right $A$-module $E$. Similarly one can define Gorenstein injective objects in the category $\widetilde{\scrR}\lRep$.
\end{ipg}

\begin{nota}
Throughout this section, denote by $\scrR$ a representation of $\calI$ on {\sf Ring} with $\scrR_i=R_i$ an associative ring for each $i\in\calI$ and $\scrR_\alpha$ a ring homomorphism for any $\alpha \in \Mor(\calI)$, and denote by
\begin{itemize}
\item $\Inj{\rMod R_j}$ the subcategory of injective right $R_j$-modules;
\item $\GI{\rMod R_j}$ the subcategory of Gorenstein injective right $R_j$-modules;
\item $\sf Inj_{\bullet}$ the family $\{ \Inj{\rMod R_j} \}_{j \in \Ob(\calI)}$ of subcategories;
\item $\sf GI_{\bullet}$ the family $\{ \GI{\rMod R_j} \}_{j \in \Ob(\calI)}$ of subcategories;
\item ${^\perp\sf GI_{\bullet}}$ the family $\{ {^\perp\GI{\rMod R_j}} \}_{j \in \Ob(\calI)}$ of subcategories;
\item $\widetilde{\scrR}$ the $\calI^{\op}$-diagram of right module categories induced by $\scrR$;
\item $\GI{\widetilde{\scrR} \lRep}$ the subcategory of Gorenstein injective objects in $\widetilde{\scrR} \lRep$;
\item $\widetilde{\scrR} \lRep_{\sf Inj_{\bullet}}$ the subcategory of objects $X\in\widetilde{\scrR} \lRep$ with $X_j \in \Inj{\rMod R_j}$ for $j \in \Ob(\calI\op)$;
\item $\widetilde{\scrR} \lRep_{^\perp \sf GI_{\bullet}}$ the subcategory of objects $X\in\widetilde{\scrR} \lRep$ with $X_j \in {^\perp \GI {\rMod R_j}}$ for $j \in \Ob(\calI\op)$.
\end{itemize}
\end{nota}

For $j \in \Ob(\calI)$, \v{S}aroch and \v{S}t'ov\'{\i}\v{c}ek \cite{SS20} constructed a hereditary Hovey triple
\[
(\rMod R_j, \, {^\perp \GI{\rMod R_j}}, \, \GI{\rMod R_j})
\]
in $\rMod R_j$ with
${^\perp \GI{\rMod R_j}} \cap \GI{\rMod R_j} = \Inj{\rMod R_j}$.
The corresponded Gorenstein injective model structure is
\begin{equation}\label{Gorenstein injective model structure}
(\Mon(\rMod R_j), \, \widetilde{\sfW}_j, \, \Epi(\GI{\rMod R_j})),
\end{equation}
where
\[\widetilde{\sfW}_j = \left \{ w_j  \left|
\begin{array}{c} w_j \text{ can be decomposed as } w_j = f_jc_j \text{ with } \\
c_j \in \Mon({^\perp \GI{\rMod R_j}}) \text{ and } f_j \in \Epi(\Inj{\rMod R_j})
\end{array} \right. \right \}.
\]

We always denote by
\begin{itemize}
\item $\overline{\Mon}(\widetilde{\scrR}) =
      \{ \sigma: N \to N' \in \widetilde{\scrR} \lRep \,|\,
      \sigma \text{ is a monomorphism}\}$,
\item $\widetilde{\sfW} =
      \{ \sigma: N \to N' \in \widetilde{\scrR} \lRep \,|\,
      \sigma_j: N_j \to N'_i \in \widetilde{\sfW}_j
      \text{ for all } j \in \Ob(\calI)\}$.
\end{itemize}

By arguments dual to the proofs of Lemmas \ref{fibrant obs} and \ref{abelian epi}, we have:

\begin{lem}\label{GI trivial obs}
Let $N$ be an object in $\widetilde{\scrR} \lRep$. Then $N \in \widetilde{\scrR} \lRep_{^\perp \sf GI_{\bullet}}$ if and only if $0 \to N \in \widetilde{\sfW}$ if and only if $N \to 0 \in \widetilde{\sfW}$.
\end{lem}

\begin{lem}\label{GI abelian mon}
Let $\sigma$ be a morphism in $\widetilde{\scrR} \lRep$. Then $\sigma \in \overline{\Mon}(\widetilde{\scrR}) \cap \widetilde{\sfW}$ if and only if $\sigma$ is a monomorphism and $0 \to \coker(\sigma) \in \overline{\Mon}(\widetilde{\scrR}) \cap \widetilde{\sfW}$.
\end{lem}

We always consider the opposite category $\calI\op$ in this section as $\widetilde{\scrR}$ is an $\calI^{\op}$-diagram. We mention that $\widetilde{\scrR}$ admits enough right adjoints, that is, each functor $\widetilde{\scrR}_{\theta}$ admits a right adjoint $\widetilde{\scrR}_{\theta}^* = \Hom[R_i]{R_j}{-}: \rMod R_i\to \rMod R_j$ for $\theta: j \to i \in \Mor(\calI\op)$; this notation will be used frequently in this section. Since we are working with the opposite category, a dual version of compatibility is required. Explicitly, let $\{\mathcal{S}_j\}_{j \in \Ob(\calI\op)}$ be a family with each $\mathcal{S}_j$ a subcategory of $\widetilde{\scrR}_j=\rMod R_j$. We say that $\{\mathcal{S}_j\}_{j \in \Ob(\calI\op)}$ is \emph{adjoint compatible with respect to} $\widetilde{\scrR}$ if $\widetilde{\scrR}_{\theta}^* (\mathcal{S}_{i}) \subseteq \mathcal{S}_{j}$ for any $\theta: j\to i \in \Mor(\calI\op)$.

\begin{lem} \label{Inj is compatible}
The family $\sf Inj_{\bullet}$ is adjoint compatible with respect to $\widetilde{\scrR}$.
\end{lem}

\begin{prf*}
This is straightforward.
\end{prf*}

The next definition will be used frequently in this section.

\begin{dfn}
A representation $\scrR$ of $\calI$ on {\sf Ring} is called {\em flat} \cite{SS2017} if $R_j$ is flat as a left $R_i$- and right $R_i$-module for any $\alpha: i \to j \in \Mor(\calI)$.
\end{dfn}

\begin{lem} \label{Ginj is compatible}
Suppose that $\scrR$ is flat. Then the family $\sf GI_{\bullet}$ is adjoint compatible with respect to $\widetilde{\scrR}$.
\end{lem}

\begin{prf*}
It follows from \lemcite[3.4]{CKFL17} that $\widetilde{\scrR}_{\alpha \op}^*(E_i) = \Hom[R_i]{R_j}{E_i}$ is Gorenstein injective in $\rMod R_j$ for any $\alpha\op: j \to i \in \Mor(\calI\op)$ and any Gorenstein injective object $E_i$ in $\rMod R_i$. The conclusion then follows.
\end{prf*}

The following result shows that under certain conditions the functor $\widetilde{\scrR}_{\theta}^*$ preserves fibrations and trivial fibrations. Therefore, $\widetilde{\scrR}_{\theta}$ preserves cofibrations and trivial cofibrations.

\begin{lem}\label{keep f and tf}
Suppose that $\scrR$ is flat. Then
\begin{prt}
\item
$\widetilde{\scrR}_{\theta}^*(\Epi(\GI{\rMod R_i})) \subseteq \Epi(\GI{\rMod R_j})$ and

\item
$\widetilde{\scrR}_{\theta}^*(\Epi(\GI{\rMod R_i}) \cap \widetilde{\sfW}_i) \subseteq \Epi(\GI{\rMod R_j}) \cap \widetilde{\sfW}_j$
\end{prt}
for $\theta: j \to i \in \Mor(\calI\op)$.
\end{lem}

\begin{prf*}
Let $f_i : A_i \to B_i \in \Epi(\GI{\rMod R_i})$. Then there exists a short exact sequence
\[
0 \to \ker(f_i) \to A_i \overset{f_i}\longrightarrow B_i \to 0
\]
with $\ker(f_i) \in \GI{\rMod R_i}$. Since $R_j$ is a flat right $R_i$-module by assumption, by \v{S}t'ov\'{\i}\v{c}ek \corcite[5.9]{JSt}, the above short exact sequence remains exact after applying the functor $\widetilde{\scrR}_{\theta}^* = \Hom[R_i]{R_j}{-}$. Thus $\widetilde{\scrR}_{\theta}^*(f_i) = \Hom[R_i]{R_j}{f_i}$ is an epimorphism. On the other hand, since the family $\sf GI_{\bullet}$ is adjoint compatible with respect to $\widetilde{\scrR}$ by Lemma \ref{Ginj is compatible}, one has $\widetilde{\scrR}_{\theta}^*(\ker(f_i)) \in \GI{\rMod R_j}$, so $\widetilde{\scrR}_{\theta}^*(f_i) \in \Epi(\GI{\rMod R_j})$.

Note that $\Epi(\GI{\rMod R_i}) \cap \widetilde{\sfW}_i = \Epi(\Inj{\rMod R_i})$ for $i \in \Ob(\calI\op)$ by \lemcite[5.8]{Ho02}. By the fact that the family $\sf Inj_{\bullet}$ is also adjoint compatible with respect to $\widetilde{\scrR}$ (see Lemma \ref{Inj is compatible}), one can prove (b) similarly.
\end{prf*}

Next, we describe a construction which is dual to the latching functor $\sfL_i$ considered in \ref{latching functor}.

\begin{ipg}\label{Matching functor}
Let $\calI$ be a direct category. Then $\calI\op$ is an inverse category. Fix $j \in \Ob(\calI\op)$, and denote by  $\calI\op_{\succ j}$ is the full subcategory of $\calI\op$ consisting of objects $i$ with $i \succ j$. For each $N \in \widetilde{\scrR} \lRep$, consider the component of the unit $N \to \iota_{\ast}\iota^{\ast}(N)$ of the adjoin pair $(\iota^*, \iota_*)$, where $\iota_{\ast}$ is the coinduction functor induced by the natural embedding functor
$\iota: \calI\op_{\succ j} \to \calI\op$; see \ref{induction-func}. Restricting to the object $j$, we get an object $\sfM_j(N) = (\iota_{\ast} \iota^{\ast} (N))_j$ as well as a natural morphism $N_j \to \sfM_j(N)$ in $\widetilde{\scrR}_j$. Then
$$\sfM_j(N) =
\lim_{\substack{\theta: \, j \to i \in \Mor(\calI\op) \\ j \neq i}}
\widetilde{\scrR}_{\theta}^*(N_j) =
\lim_{\theta \in \calP_j(j, \bullet)} \widetilde{\scrR}_{\theta}^*(N_\bullet),$$
where $\calP_j=\Mor(\calI\op)\backslash{\sf End}_{\calI\op}(j)$. It is easy to see that $\sfM_j$ is precisely the matching functor described in \cite{HR08}, and the morphism $N_j \to \sfM_j(N)$ is actually the morphism $\psi_j^N: N_j \to \lim_{\theta \in \calP_j(j, \bullet)} \widetilde{\scrR}_{\theta}^*(N_\bullet)$ given in \cite[(2.17.1)]{DLLY}.
\end{ipg}

\begin{dfn} \label{t and psi}
Given a direct category $\calI$ (in this case $\calI\op$ is an inverse category) and a family $\calY = \{ \calY_j \}_{j \in \Ob(\calI\op)}$ with each $\calY_j$ a full subcategory of $\widetilde{\scrR}_j$, define a subcategory of $\widetilde{\scrR}\lRep$:
\[
\Psi(\calY)= \{ N \in \widetilde{\scrR}\lRep \mid \psi_j^N \text{ is an epimorphism and } \ker(\psi_j^N) \in {\calY}_j \text{ for each } j \in \Ob(\calI\op) \}.
\]
\end{dfn}

For any morphism $\sigma: N \to K \in \widetilde{\scrR} \lRep$,
considering the following commutative diagram in which the inner square is a pullback, we get a natural morphism $\varrho_j$:
\begin{equation} \label{Match pullback}
\xymatrix{
  N_j \ar@/_/[ddr]_{\sigma_j} \ar@{.>}[dr]|-{\varrho_j} \ar@/^/[drr]^{\psi_j^{N}}\\
& K_j \times_{\sfM_j(K)} \sfM_j(N) \ar[r]_-{\varsigma_j} \ar[d]^-{\kappa_j} & \sfM_j(N) \ar[d]^{\sfM_j(\sigma)}\\
& K_j \ar[r]_-{\psi_j^{K}} & \sfM_j(K).  }
\end{equation}
Set $\overline{\PBEpi}({\sf GI_{\bullet}}) =
      \left \{ \sigma: N \to K \in \widetilde{\scrR} \lRep \, \left|
      \begin{array}{c}
      \varrho_j: N_j \to K_j \times_{\sfM_j(K)} \sfM_j(N) \in \Epi(\GI{\rMod R_j})\\
      \text{ for all } j \in \Ob(\calI\op) \end{array} \right. \right \}$.

With help of Lemma \ref{keep f and tf}, we obtain the following result.

\begin{prp}\label{prp5.12}
Suppose that $\calI$ is direct, and $\scrR$ is flat. Then the triplet
$$(\overline{\Mon}(\widetilde{\scrR}), \, \widetilde{\sfW}, \, \overline{\PBEpi}({\sf GI_{\bullet}}))$$
is a model structure on $\widetilde{\scrR} \lRep$.
\end{prp}
\begin{prf*}
We mention that $(\Mon(\rMod R_j), \, \widetilde{\sfW}_j, \, \Epi(\GI{\rMod R_j}))$ is an abelian model structure on $\rMod R_j$ for each $j \in \Ob(\calI)$; see (\ref{Gorenstein injective model structure}). By Lemma \ref{keep f and tf}, one gets that $\widetilde{\scrR}_{\theta}^*$ preserves fibrations and trivial fibrations for each $\theta: j \to i \in \Mor(\calI\op)$, and so $\widetilde{\scrR}_{\theta}$ preserves cofibrations and trivial cofibrations. It is easy to see that the triplet
$$(\overline{\Mon}(\widetilde{\scrR}), \, \widetilde{\sfW}, \, \overline{\PBEpi}({\sf GI_{\bullet}}))$$
is the $f$-structure in the sense of H\"{u}ttemann and R\"{ondigs}; see \cite[Definition 3.3.3]{HR08}. Note that $\widetilde{\scrR}$ is an $\calI\op$-diagram, while $\calI\op$ is inverse. Then it follows from \thmcite[3.3.5(1)]{HR08} that the above triplet is a model structure on $\widetilde{\scrR} \lRep$.
\end{prf*}

Next, we show that the model structure on $\widetilde{\scrR} \lRep$ in Proposition \ref{prp5.12} is abelian whenever $\calI$ is a left rooted quiver (viewed as a category in a natural way). We mention that in this case
\[
\sfM_j(N)=\lim_{\theta \in \calP_j(j, \bullet)} \widetilde{\scrR}_{\theta}^*(N_\bullet) = \prod_{\theta \in \calI\op(j, \bullet)} \widetilde{\scrR}_{\theta}^*(N_\bullet),
\]
where $\calI\op(j, \bullet)$ denotes the set of all arrows in $\calI\op$ with $j$ the source, and the last equality follows from a dual of \cite[Example 2.3]{DLLY}. We also mention that the morphism $\psi_j^N$ given in \ref{Matching functor} is from $N_j$ to $\prod_{\theta \in \calI\op(j, \bullet)} \widetilde{\scrR}_{\theta}^*(N_\bullet)$.
The main obstacle forcing us to work with left rooted quiver rather than arbitrary direct categories is that the limit functor in general is not exact, but this exactness is essential for us to construct Gorenstein model structures. For left rooted quivers, the limit functor coincides with the product functor, which is exact.

\begin{lem}\label{fibration epi}
Suppose that $\calI$ is a left rooted quiver, and $\scrR$ is flat. Then a morphism $\sigma$ in $\overline{\PBEpi}({\sf GI_{\bullet}})$ is an epimorphism
with $\ker(\sigma_j) \in \GI{\rMod R_j}$ for each $j \in \Ob(\calI\op)$.
\end{lem}

\begin{prf*}
Let $\{V_{\chi}\}_{\chi\,{\rm ordinal}}$ be the transfinite sequence of subsets of $\Ob(\calI)$ defined in \ref{construct}. Then there exists a certain ordinal $\zeta$ such that $\Ob(\calI) = V_{\zeta}$. Take $\sigma: N \to K \in \overline{\PBEpi}({\sf GI_{\bullet}})$. We use the transfinite induction to show the following conclusion: $\sigma_j: N_j \to K_j$ is an epimorphism with $\ker(\sigma_j) \in \GI{\rMod R_j}$ for all $j \in V_{\chi}$ and all ordinals $\chi \leqslant \zeta$.

Note that $\varrho_j$ is always an epimorphism with $\ker(\varrho_j)\in \GI{\rMod R_j}$ for all $j \in \Ob(\calI)$. Consider the commutative diagram (\ref{Match pullback}). If $j \in V_1$, then it is minimal with respect to the partial order $\preccurlyeq$ in $\calI$, or equivalently, there exists no arrow in $\calI\op$ such that $j$ is its source. It follows that $\sfM_j (N) = 0 = \sfM_j (K)$, $\sfM_j (\sigma) = 0$ and $\kappa_j = \id_{K_j}$. Consequently, $\sigma_j = \varrho_j$ which satisfies the conclusion.

For $\chi > 1$, we have two cases:

(1) $\chi$ is a successor ordinal. Take $j \in V_{\chi}$ and let $i \in \Ob(\calI)$ with $i \prec j$ in $\Mor(\calI)$, that is, $i \neq j$ and there exists an arrow  $\theta : j \to i$ in $\calI\op$. Then we have $i \in V_{\chi-1}$ by \cite[Remark3.3]{DLLY}.  By the induction hypothesis, all $\sigma_i$ are epimorphisms with $\ker(\sigma_i) \in \GI{\rMod R_i}$.  Since $R_j$ is a flat right $R_i$-module, it follows from \corcite[5.9]{JSt} that $\widetilde{\scrR}_{\theta}^*(\sigma_i) = \Hom[R_i]{R_j}{\sigma_i}$ is an epimorphism. Consequently,
$\sfM_j(\sigma)=\prod_{\theta \in \calI\op(j, \bullet)} \widetilde{\scrR}_{\theta}^*(\sigma_\bullet)$
is an epimorphism, so $\kappa_j$ is an epimorphism, and so is $\sigma_j = \kappa_j \varrho_j$.
On the other hand, note that the family $\sf GI_{\bullet}$ is adjoint compatible with respect to $\widetilde{\scrR}$ by Lemma \ref{Ginj is compatible}, so $\widetilde{\scrR}_{\theta}^*(\ker(\sigma_i)) \in \GI{\rMod R_j}$. Thus one has
\[
\ker(\sfM_j(\sigma)) = \ker(\prod_{\theta \in \calI\op(j, \bullet)} \widetilde{\scrR}_{\theta}^*(\sigma_\bullet)) \cong \prod_{\theta \in \calI\op(j, \bullet)} \ker(\widetilde{\scrR}_{\theta}^*(\sigma_\bullet))\cong \prod_{\theta \in \calI\op(j, \bullet)} \widetilde{\scrR}_{\theta}^*(\ker\sigma_\bullet)\in \GI{\rMod R_j},
\]
which implies that $\ker(\kappa_j) \in \GI{\rMod R_j}$ as $\ker(\kappa_j) \cong \ker(\sfM_j(\sigma))$.  Consider now the following commutative diagram with exact rows and columns
\[
\xymatrix@C=20pt@R=15pt{
& & \ker(\varrho_j) \ar@{>->}[d] \\
0 \ar[r] & \ker(\sigma_j) \ar[r] \ar[d] & N_j \ar[r]^{\sigma_j} \ar@{>>}[d]^{\varrho_j} & K_j \ar[r] \ar@{=}[d] & 0 \\
0 \ar[r] & \ker(\kappa_j) \ar[r] & K_j \times_{\sfM_j(K)} \sfM_j(N)  \ar[r]^-{\kappa_j} & K_j \ar[r] & 0. \\
}
\]
By the snake lemma, one gets a short exact sequence $0 \to \ker(\varrho_j) \to \ker(\sigma_j) \to \ker(\kappa_j) \to 0$ in $\widetilde{\scrR}_j$. Since both $\ker(\varrho_j)$ and $\ker(\kappa_j)$ are in $\GI{\rMod R_j}$, so is $\ker(\sigma_j)$. Thus $\sigma_j$ also satisfies the conclusion.

(2) If $\chi \leqslant \zeta$ is a limit ordinal, then the assertion is clearly true for $\chi$ because in this case $V_{\chi} = \cup_{\mu < \chi} V_{\mu}$ and the conclusion holds for all ordinals $\mu < \chi$.

Finally, the conclusion follows by taking $\chi = \zeta$.
\end{prf*}

\begin{prp}\label{Rondigs f}
Suppose that $\calI$ is a left rooted quiver, and $\scrR$ is flat. Then
\[
(\overline{\Mon}(\widetilde{\scrR}), \, \widetilde{\sfW}, \, \overline{\PBEpi}({\sf GI_{\bullet}}))
\]
is an abelian model structure on $\widetilde{\scrR} \lRep$.
\end{prp}
\begin{prf*}
It follows from Proposition \ref{prp5.12} that $(\overline{\Mon}(\widetilde{\scrR}), \, \widetilde{\sfW}, \, \overline{\PBEpi}({\sf GI_{\bullet}}))$ is a model structure on $\widetilde{\scrR} \lRep$. It is clear that cofibrations are monomorphisms and fibrations are epimorphisms by Lemma \ref{fibration epi}. To prove the above model structure is abelian, we have to prove that cofibrations coincide with monomorphisms with cofibrant cokernels, and trivial cofibrations coincide with monomorphisms with trivial cofibrant cokernels; see Remark \ref{df of ams is self-dual}. However, the first statement holds clearly, and second one follows from Lemma \ref{GI abelian mon}.
\end{prf*}

We will show later that the abelian model structure in Proposition \ref{Rondigs f} is precisely the Gorenstein injective model structure on $\widetilde{\scrR} \lRep$.

\begin{lem}\label{fibrant obs2}
Suppose that $\calI$ is a left rooted quiver.
Then an object $N \in \Psi(\sf GI_{\bullet})$ if and only if
$N \to 0 \in \overline{\PBEpi}({\sf GI_{\bullet}})$.
\end{lem}

\begin{prf*}
Consider the commutative diagram (\ref{Match pullback}) with $K = 0$. Then $\varsigma_j = \id_{\sfM_j(N)}$, and so $\psi_j^{N} = \varrho_j$. Thus $N \in \Psi(\sf GI_{\bullet})$ if and only if $\varrho_j$ is an epimorphism
with $\ker (\varrho_j) \in \GI{\rMod R_j}$ for each $j \in \Ob(\calI\op)$. But the later statement is equivalent to saying that $N \to 0 \in \overline{\PBEpi}({\sf GI_{\bullet}})$.
\end{prf*}

\begin{thm}\label{ht induce ht GI 1}
Suppose that $\calI$ is a left rooted quiver, and $\scrR$ is flat. Then
\[
(\widetilde{\scrR} \lRep, \, \widetilde{\scrR} \lRep_{^\perp \sf GI_{\bullet}}, \, \Psi(\sf GI_{\bullet}))
\]
forms a hereditary Hovey triple in $\widetilde{\scrR} \lRep$.
\end{thm}

\begin{prf*}
Since $(\overline{\Mon}(\widetilde{\scrR}), \, \widetilde{\sfW}, \, \overline{\PBEpi}({\sf GI_{\bullet}}))$ forms an abelian model structure on $\widetilde{\scrR} \lRep$ by Proposition \ref{Rondigs f}, we only need to use Theorem \ref{Hovey cor} to recognize the associated subcategories of the cofibrant, trivial and fibrant objects. We have:
\begin{itemize}
\item the subcategory of cofibrant objects is obviously $\widetilde{\scrR} \lRep$;

\item the subcategory of trivial objects is $\widetilde{\scrR} \lRep_{^\perp \sf GI_{\bullet}}$ by Lemma \ref{GI trivial obs};

\item the subcategory of fibrant objects is $\Psi(\sf GI_{\bullet})$ by Lemma \ref{fibrant obs2}.
\end{itemize}
Thus $(\widetilde{\scrR} \lRep, \, \widetilde{\scrR} \lRep_{^\perp \sf GI_{\bullet}}, \, \Psi(\sf GI_{\bullet}))$ forms a Hovey triple in $\widetilde{\scrR} \lRep$.

To show the above Hovey triple is hereditary, we have to prove the complete cotorsion pairs
\[
(\widetilde{\scrR} \lRep_{^\perp \sf GI_{\bullet}}, \, \Psi(\sf GI_{\bullet})) \text{ and } (\widetilde{\scrR} \lRep, \, \widetilde{\scrR} \lRep_{^\perp \sf GI_{\bullet}} \cap \Psi(\sf GI_{\bullet}))
 \]
in $\widetilde{\scrR} \lRep$ are hereditary. By Lemma \ref{test hereditary}, it is enough to show that both the above two cotorsion pairs are resolving.
But this is obvious.
\end{prf*}

\begin{rmk}
Let $\calI$ be a left rooted quiver and $\scrD$ a $\calI\op$-diagram admitting enough right adjoints. If the right adjoint $\scrD_\theta^*$ of $\scrD_\theta$ is exact for any arrow $\theta$ in $\calI\op$, then by a dual way for obtaining Theorem \ref{ht induce ht 1}, one can show that any family
\[
\{(\calQ_j, \, \calW_j, \, \calR_j)\}_{j \in \Ob(\calI\op)}
\]
of Hovey triples in $\scrD_j$'s such that both $\calR$ and $\widetilde{\calR}$ are adjoint compatible with respect to $\scrD$ induces a Hovey triple
\[
(\scrD \lRep_{\calQ}, \, \scrD \lRep_{\calW}, \, \Psi(\calR))
\]
in $\scrD \lRep$. However, Theorem \ref{ht induce ht GI 1} is not a special case of the above assertion, though it is proved in a similar way. Indeed, even if $\scrR$ is flat, the right adjoint $\widetilde{\scrR}_\theta^*$ of $\widetilde{\scrR}_\theta$ is not exact in general, that is, $\widetilde{\scrR}_\theta^* = \Hom[R_i]{R_j}{-}$ might not be exact for every arrow $\theta: j \to i$ in $\calI\op$. The key point we used to obtain Theorem \ref{ht induce ht GI 1} is that $\Ext[R_i]{1}{R_j}{E_i} = 0$ for each Gorenstein injective right $R_i$-module $E_i$; see \corcite[5.9]{JSt}.
\end{rmk}

As an immediate consequence of Theorem \ref{ht induce ht GI 1},
we obtain the following result.

\begin{cor} \label{Ginj cp over left rooted}
Suppose that $\calI$ is a left rooted quiver and $\scrR$ is flat. Then $(\widetilde{\scrR} \lRep_{^\perp \sf GI_{\bullet}}, \, \Psi(\sf GI_{\bullet}))$ is a complete and hereditary cotorsion pair in $\widetilde{\scrR} \lRep$.
\end{cor}

Our next task is to show that under some assumptions the subcategories $\GI{\widetilde{\scrR} \lRep}$ and $\Psi(\sf GI_{\bullet})$ of ${\widetilde{\scrR} \lRep}$ coincide. Before proving this result, we need to finish a few preparatory works.

\begin{lem}\label{X(i) in X inj}
Suppose that $\calI$ is a left rooted quiver. Then one has $\Psi(\sf Inj_{\bullet}) \subseteq \widetilde{\scrR} \lRep_{\sf Inj_{\bullet}}$.
\end{lem}

\begin{prf*}
Note that the family $\sf Inj_{\bullet}$ is adjoint compatible with respect to $\widetilde{\scrR}$ by Lemma \ref{Inj is compatible}, and $\Inj{\rMod R_j}$ is closed under
extensions and small products for each $j \in \Ob(\calI)$.
By a dual way for proving Lemma \ref{X(i) in X},
the desired result follows.
\end{prf*}

In the next result we collect some elementary properties of Gorenstein injective objects in $\widetilde{\scrR} \lRep$.

\begin{lem} \label{on side for g injectives}
Suppose that $\calI$ is a left rooted quiver, and $\scrR$ is flat. Let $N$ be a Gorenstein injective object in $\widetilde{\scrR} \lRep$. Then for any $j \in \Ob(\calI\op)$, ${\psi}_j^{N}: N_j \to \prod_{\theta \in \calI\op(j,\bullet)} \Hom[R_{t(\theta)}]{R_j}{N_{t(\theta)}}$ is an epimorphism with $\ker({\psi}_j^{N})$ Gorenstein injective in $\rMod R_j$. That is, there is a containment $\GI{\widetilde{\scrR} \lRep} \subseteq \Psi(\sf GI_{\bullet})$.
\end{lem}

\begin{prf*}
Fix $j \in \Ob(\calI\op)$. We define a functor $\K_j: \widetilde{\scrR}\lRep\to\rMod R_j$ sending a representation $M\in\widetilde{\scrR}\lRep$ to $\ker({\psi}_j^{M})$; see \cite[Corollary 2.19]{DLLY}. Since $N$ is a Gorenstein injective object in $\widetilde{\scrR} \lRep$, there exists an exact sequence
\[
\mathbb{I}: \, \cdots \to I^{-1} \to I^0 \to I^1 \to \cdots
\]
of injective objects in $\widetilde{\scrR}\lRep$ such that $N \cong \ker{(I^0 \to I^1)}$ and the sequence $\mathbb{I}$ remains exact after applying the functor $\Hom[\widetilde{\scrR}\lRep]{E}{-}$ for every injective object $E$ in $\widetilde{\scrR}\lRep$.
For each $i \in \mathbb{Z}$, since $I^i \in \Psi(\sf Inj_{\bullet})$ by \cite[Corollary 3.20]{DLLY}, there exists a short exact sequence
\[
0 \to \K_j(I^i) \to I^i_j \to \prod_{\theta \in \calI\op(j, \bullet)}\Hom[R_{t(\theta)}]{R_j}{I^i_{t(\theta)}} \to 0
\]
in $\rMod R_j$ with $\K_j(I^i)$ injective. For any arrow $\theta \in \calI\op(j, \bullet)$,
we have an exact sequence
\[
\mathbb{I}_{t(\theta)}= \quad \cdots \to I^{-1}_{t(\theta)} \to I^0_{t(\theta)} \to I^{1}_{t(\theta)} \to \cdots
\]
in $\rMod R_{t(\theta)}$. By Lemma \ref{X(i) in X inj}, all $I^i_{t(\theta)}$ are injective in $\rMod R_{t(\theta)}$. It follows from \corcite[5.9]{JSt} that the sequence $\Hom[R_{t(\theta)}] {R_j} {\mathbb{I}_{t(\theta)}}$ is exact as $R_j$ is flat in $\rMod R_{t(\theta)}$, so the sequence $\prod_{\theta \in \calI\op(j,\bullet)} \Hom[R_{t(\theta)}] {R_j} {\mathbb{I}_{t(\theta)}}$ is exact. Consequently, we get the commutative diagram
\[
\xymatrix@C=20pt@R=15pt{
 & \vdots \ar[d] & \vdots \ar[d] & \vdots \ar[d] \\
 0 \ar[r] & \K_j(I^{-1}) \ar[d] \ar[r] & I^{-1}_j \ar[d] \ar[r] &
\prod_{\theta \in \calI\op(j,\bullet)} \Hom[R_{t(\theta)}] {R_j} {I^{-1}_{t(\theta)}}
\ar[d] \ar[r] & 0 \\
 0 \ar[r] & \K_j(I^0) \ar[d] \ar[r] & I^0_j \ar[d] \ar[r] &
 \prod_{\theta \in \calI\op(j,\bullet)} \Hom[R_{t(\theta)}] {R_j} {I^0_{t(\theta)}}
 \ar[d] \ar[r] & 0 \\
 0 \ar[r] & \K_j(I^{1}) \ar[d] \ar[r] & I^{1}_j \ar[d] \ar[r] &
 \prod_{\theta \in \calI\op(j,\bullet)} \Hom[R_{t(\theta)}] {R_j} {I^{1}_{t(\theta)}}
 \ar[d] \ar[r] & 0 \\
 & \vdots & \vdots & \vdots
}
\]
with exact rows and columns, which induces the short exact sequence
\[
0 \to \K_j(N) \to N_j \overset{{\psi}_j^{N}} \longrightarrow
\prod_{\theta \in \calI\op(j,\bullet)} \Hom[R_{t(\theta)}] {R_j} {N_{t(\theta)}} \to 0.
\]
Therefore, to complete the proof, it remains to show that $\K_j(N)$ is Gorenstein injective in $\rMod R_j$, that is, the sequence $\Hom[R_j] {E_j} {\K_j(\mathbb{I})}$ is exact for any injective object $E_j$ in $\rMod R_j$.

Indeed, we have
\begin{align*}
\Hom[\widetilde{\scrR} \lRep] {\cofre_j(E_j)} {\mathbb{I}}
& \cong
\Hom[\widetilde{\scrR} \lRep] {\cofre_j(E_j)}
{\prod_{h \in \Ob(\calI\op)}\cofre_h(\K_h(\mathbb{I}))}\\
& \cong
\prod_{h \in \Ob(\calI\op)} \Hom[\widetilde{\scrR} \lRep] {\cofre_j(E_j)} {\cofre_h(\K_h(\mathbb{I}))}\\
& \cong
\prod_{h \in \Ob(\calI\op)} \Hom[R_h] {\eva^h(\cofre_j(E_j))} {\K_h(\mathbb{I})},
\end{align*}
where the first isomorphism holds by \cite[Remark 3.22]{DLLY} and the third isomorphism holds as $(\eva^h, \cofre_h)$ is an adjoint pair; see \ref{free-func}. Since $\cofre_j(E_j)$ is injective in $\widetilde{\scrR} \lRep$ as the functor $\eva^j$ is exact, the sequence $\Hom[\widetilde{\scrR} \lRep] {\cofre_j(E_j)} {\mathbb{I}}$ is exact, and so $\Hom[R_j] {\eva^j(\cofre_j(E_j))} {\K_j(\mathbb{I})}$ is also exact by the above isomorphisms. The conclusion then follows by observing that
$\eva^j (\cofre_j(E_j)) = \widetilde{\scrR}_{e_j}(E_j) = E_j \otimes_{R_j} R_j \cong E_j$.
\end{prf*}

Let $\calY$ be a subcategory of an abelian category $\calA$.
Recall from Auslander and Buchweitz \cite{MAsROB89} that
a subcategory $\calV$ of $\calY$ is called a \emph{generator} for $\calY$ if for each object $Y \in \calY$, there exists a short exact sequence $0 \to Y' \to V \to Y \to 0$ with $V \in \calV$ and $Y' \in \calY$.

\begin{lem} \label{in cogenertor}
Suppose that $\calI$ is a left rooted quiver, and $\scrR$ is flat. Then $\Psi(\sf Inj_{\bullet})$ is a generator for $\Psi(\sf GI_{\bullet})$.
\end{lem}

\begin{prf*}
Take $N \in \Psi(\sf GI_{\bullet})$. Since $(\widetilde{\scrR} \lRep_{^\perp \sf GI_{\bullet}}, \, \Psi(\sf GI_{\bullet}))$ is a complete cotorsion pair in $\widetilde{\scrR} \lRep$ by Corollary \ref{Ginj cp over left rooted}, there exists a short exact sequence
$0 \to N' \to E \to N \to 0$ in $\widetilde{\scrR} \lRep$ with $E \in \widetilde{\scrR} \lRep_{^\perp \sf GI_{\bullet}} \cap \Psi(\sf GI_{\bullet})$ and $N' \in \Psi(\sf GI_{\bullet})$. By \cite[Corollary 3.20]{DLLY} and Theorem \ref{ht induce ht GI 1},
\[
\widetilde{\scrR} \lRep_{^\perp \sf GI_{\bullet}} \cap \Psi(\sf GI_{\bullet}) = \Inj{\widetilde{\scrR} \lRep} = \Psi(\sf Inj_{\bullet}),
\]
so $E \in \Psi({\sf Inj_{\bullet}})$. This finishes the proof.
\end{prf*}

Now we are ready to give a characterization of Gorenstein injective objects in $\widetilde{\scrR} \lRep$, which generalizes a result by Eshraghi, Hafezi and Salarian; see \cite[Theorem 3.5.1]{EHS13}.

\begin{thm} \label{Gorenstein in r module}
Suppose that $\calI$ is a left rooted quiver, and $\scrR$ is flat. Then $\GI{\widetilde{\scrR} \lRep} = \Psi(\sf GI_{\bullet})$.
\end{thm}

\begin{prf*}
The inclusion $\GI{\widetilde{\scrR} \lRep} \subseteq \Psi(\sf GI_{\bullet})$ holds by Lemma \ref{on side for g injectives}. For the other inclusion, we take $N \in \Psi(\sf GI_{\bullet})$ and show that $N \in \GI{\widetilde{\scrR} \lRep}$.

Since $\Inj{\widetilde{\scrR} \lRep} = \Psi(\sf Inj_{\bullet})$ is a generator for $\Psi(\sf GI_{\bullet})$ by \cite[Corollary 3.20]{DLLY} and Lemma \ref{in cogenertor}, there exists a short exact sequence
$0 \to N^{-1} \to I^{-1} \to  N \to 0$ in $\widetilde{\scrR} \lRep$ with $I^{-1} \in \Inj{\widetilde{\scrR} \lRep}$ and $N^{-1} \in \Psi(\sf GI_{\bullet})$. For any injective object $E$ in $\widetilde{\scrR} \lRep$, by Lemma \ref{X(i) in X inj}, $E \in \widetilde{\scrR} \lRep_{\sf Inj_{\bullet}} \subseteq \widetilde{\scrR} \lRep_{^\perp \sf GI_{\bullet}}$. Thus by Corollary \ref{Ginj cp over left rooted}, one has $\Ext[{\widetilde{\scrR} \lRep}]{1}{E}{N^{-1}} = 0$, and so the above short exact sequence remains exact after applying the functor $\Hom[{\widetilde{\scrR} \lRep}]{E}{-}$. Continuing this process for $N^{-1}$, eventually one gets an exact sequence
\[
\cdots \to I^{-2} \to I^{-1} \to N \to 0 \quad \quad \quad (\dag)
\]
in $\widetilde{\scrR} \lRep$ with $I^i \in \Inj{\widetilde{\scrR} \lRep}$ for negative integers $i$ such that the sequence $(\dag)$ remains exact after applying the functor $\Hom[{\widetilde{\scrR} \lRep}]{E}{-}$.

On the other hand, since $\widetilde{\scrR} \lRep$ has enough injectives, we get an exact sequence
\[
0 \to N \to I^0 \to I^{1} \to \cdots \quad \quad \quad (\ddag)
\]
in $\widetilde{\scrR} \lRep$ with $I^i \in \Inj{\widetilde{\scrR} \lRep}$ for all integers $i \geqslant 0$. By Corollary \ref{Ginj cp over left rooted} again, $\Psi(\sf GI_{\bullet})$ is closed under taking cokernels of monomorphisms, so all cokernels of the sequence $(\ddagger)$ belong to $\Psi(\sf GI_{\bullet})$ by noting that each $I^i \in \Psi(\sf Inj_{\bullet})\subseteq \Psi(\sf GI_{\bullet})$. Using the same argument as above, we conclude that the sequence $(\ddag)$ remains exact after applying the functor $\Hom[{\widetilde{\scrR} \lRep}]{E}{-}$.

Assembling the exact sequences $(\dagger)$ and $(\ddagger)$, we deduce that $N \in \GI{\widetilde{\scrR} \lRep}$, as desired.
\end{prf*}

As an immediate consequence of
Corollary \ref{Ginj cp over left rooted} and Theorem \ref{Gorenstein in r module},
we have:

\begin{cor} \label{cp for ginj}
Suppose that $\calI$ is a left rooted quiver, and $\scrR$ is flat. Then $(\widetilde{\scrR} \lRep_{^\perp{\sf GI_{\bullet}}}, \, \GI{\widetilde{\scrR} \lRep})$ is a complete and hereditary cotorsion pair in $\widetilde{\scrR} \lRep$.
\end{cor}

The following result gives a Gorenstein injective model structure on $\widetilde{\scrR} \lRep$.

\begin{cor} \label{ginj model structor}
Suppose that $\calI$ is a left rooted quiver, and $\scrR$ is flat. Then there is a hereditary Hovey triple
$(\widetilde{\scrR} \lRep, \, {^\perp \GI{\widetilde{\scrR} \lRep}}, \, \GI{\widetilde{\scrR} \lRep})$ in $\widetilde{\scrR} \lRep$.
\end{cor}

\begin{prf*}
By Propositin \ref{ht induce ht GI 1}, there exists a hereditary Hovey triple
$(\widetilde{\scrR} \lRep, \, (\widetilde{\scrR} \lRep)_{^\perp \sf GI_{\bullet}}, \, \Psi(\sf GI_{\bullet}))$
in $\widetilde{\scrR} \lRep$. The conclusion then follows from Theorem \ref{Gorenstein in r module} and Corollary \ref{cp for ginj}.
\end{prf*}

\section{A characterization of flat objects in $\overline{\scrR} \lRep$}\label{flat}
\noindent
We now turn our attention to the category $\overline{\scrR} \lRep$, where $\overline{\scrR}$ is the $\calI$-diagram of left module categories induced by a representation $\scrR$ of $\calI$ on {\sf Ring}, with $\overline{\scrR}_i=R_i\lMod$ for $i \in \Ob(\calI)$ and $\overline{\scrR}_{\alpha}=R_j\otimes_{R_i}-: R_i\lMod \to R_j\lMod$ for $\alpha: i \to j \in \Mor(\calI)$; see Example \ref{bar R}.Note that $\overline{\scrR} \lRep$ is a locally finitely presented Grothendieck category admitting enough projectives, by Theorem \ref{Rep is grothendieck}. To characterize Gorenstein flat objects in the category $\overline{\scrR} \lRep$ and establish Gorenstein flat model structures, we provide a description of the categorical flat objects in $\overline{\scrR} \lRep$ in this section based on the following definition of the categorical tensor products due to Oberst and Rohrl \cite{Oberst70}.

\begin{ipg}\label{categorical tensor product}
Given an object $M$ in $\overline{\scrR} \lRep$ and a $\bbZ$-module $G$, define $\Hom{M}{G}$ as follows:
\begin{itemize}
\item for $i \in \Ob(\calI\op)$, set $\Hom{M}{G}_i$ to be $\Hom[\bbZ]{M_i}{G} \in \rMod R_i$;

\item for $\alpha\op: j \to i \in \Mor(\calI\op)$, there exists a morphism
\[
M_\alpha :\overline{\scrR}_\alpha(M_i) = R_j \otimes_{R_i} M _i \to M_j
\]
in $\overline{\scrR}_j=R_j\lMod$, which yields a morphism
\[\quad\quad\quad
\Hom[\bbZ]{M_\alpha}{G}: \Hom[\bbZ]{M_j}{G} \longrightarrow \Hom[\bbZ]{R_j \otimes_{R_i} M _i}{G}\overset{\cong} \longrightarrow\Hom[R_i]{R_j}{\Hom[\bbZ]{M_i}{G}}
\]
in $\widetilde{\scrR}_j=\rMod R_j$. Set $\Hom{M}{G}_{\alpha\op}$ to be the adjoint morphism of $\Hom[\bbZ]{M_\alpha}{G}$ with respect to the adjoint pair $(-\otimes_{R_j}R_j, \Hom[R_i]{R_j}{-})$.
\end{itemize}

It is routine to check that $\Hom{M}{G}$ is in $\widetilde{\scrR}\lRep$, and $\Hom{M}{-}$ is a functor from $\bbZ\lMod $ to $\widetilde{\scrR} \lRep$; it is left exact and preserves products. Consequently, the functor $\Hom{M}{-}$ has a left adjoint, denoted by $\tp[R]{M}{-}: \widetilde{\scrR} \lRep \to \bbZ \lMod$, and is called the \textit{tensor product functor}. Given an object $N\in\widetilde{\scrR}\lRep$, the tensor product functor $\tp[R]{-}{N}: \overline{\scrR}\lRep \to \bbZ \lMod$ can be defined similarly.
\end{ipg}

The next result is obtained immediately.

\begin{lem} \label{adjoint-iso for tensor product}
Let $M$ be an object in $\overline{\scrR} \lRep$ and $N$ an object in $\widetilde{\scrR} \lRep$. Then for any $\bbZ$-module $G$, there exists a natural isomorphism
\[
\Hom[\bbZ]{\tp[R]{M}{N}}{G} \cong \Hom[\widetilde{\scrR} \lRep]{N}{\Hom{M}{G}}.
\]
\end{lem}

In the following we give the definition of categorical flat objects in $\overline{\scrR} \lRep$.

\begin{dfn} \label{definition of flat}
An object $F$ in ${\overline{\scrR} \lRep}$ is called \emph{flat} if the functor $\tp[R]{F}{-}$ is exact. The subcategory of all flat objects in $\overline{\scrR}\lRep$ is denoted by $\F{\overline{\scrR}\lRep}$.
\end{dfn}

An equivalent characterization of categorical flat object is:

\begin{lem} \label{equivalent for flat}
An object $F$ in $\overline{\scrR} \lRep$ is flat if and only if $F^+ = \Hom[\bbZ]{F}{\bbQ/\bbZ}$ is injective in $\widetilde{\scrR} \lRep$.
\end{lem}

\begin{prf*}
Given a short exact sequence $0 \to N \to N' \to N'' \to 0$ in $\widetilde{\scrR} \lRep$, we consider the following commutative diagram
\[
\xymatrix@C=15pt@R=25pt{
0 \ar[r] & \Hom[\widetilde{\scrR} \lRep]{N}{F^+} \ar[r] \ar[d]^{\cong} &
\Hom[\widetilde{\scrR} \lRep]{N'}{F^+} \ar[r] \ar[d]^{\cong} &
\Hom[\widetilde{\scrR} \lRep]{N''}{F^+} \ar[r] \ar[d]^{\cong} & 0 \\
0 \ar[r] & (\tp[R]{F}{N})^+ \ar[r] & (\tp[R]{F}{N'})^+ \ar[r] & (\tp[R]{F}{N''})^+ \ar[r] & 0
}
\]
of abelian groups, where the columns are isomorphisms by Lemma \ref{adjoint-iso for tensor product}. Then the first row is exact if and only if so is the second one, or equivalently, $F$ is flat if and only if $F^+$ is injective.
\end{prf*}

Let ${\sf Flat}_{\bullet} = \{ \F {R_i \lMod} \}_{i \in \Ob(\calI)}$ be the family of subcategories of flat left $R_i$-modules. For a object $j\in\Ob(\calI)$, we mention that $\calP_j=\Mor(\calI)\backslash\AU j$, where $\AU j$ is the set of endomorphisms on $j$.

\begin{lem} \label{corresponding relation}
Suppose that $\calI$ is a direct category. Then an object $M$ in $\overline{\scrR} \lRep$ is contained in the subcategory $\Phi(\sf Flat_{\bullet})$ of $\overline{\scrR} \lRep$ if and only if $M^+$ is contained in the subcategory $\Psi(\sf Inj_{\bullet})$
of $\widetilde{\scrR} \lRep$.
\end{lem}
\begin{prf*}
The statement that $M \in \Phi(\sf Flat_{\bullet})$ is equivalent to the second statement:
for all $j \in \Ob(\calI)$, there exists a short exact
\[
0 \to \colim_{\sigma \in \calP_j(\bullet,j)}(R_j \otimes_{R_{s(\sigma)}} M_{s(\sigma)}) \overset{{\varphi}_j^{M}} \longrightarrow M_j \to \coker({\varphi}_j^{M}) \to 0
\]
in $R_j \lMod$ with $\coker({\varphi}_j^{M})$ flat; see Definition \ref{s and phi}. This turns out to be equivalent to the third statement: for $j \in \Ob(\calI)$, there exists a short exact sequence
\[
0 \to (\coker({\varphi}_j^{M}))^+ \to (M_j)^+ \overset{({\varphi}_j^M)^+} \longrightarrow (\colim_{\sigma \in \calP_j(\bullet,j)} (R_j \otimes_{R_{s(\alpha)}} M_{s(\alpha)}))^+ \to 0
\]
in $\rMod R_j$ with $(\coker({\varphi}_j^{M}))^+$ injective. However, since
\begin{align*}
(\colim_{\sigma \in \calP_j(\bullet,j)} (R_j \otimes_{R_{s(\sigma)}} M_{s(\sigma)}))^+
& = \Hom[\bbZ]{\colim_{\sigma\op \in \calP_j(j,\bullet)}(R_j \otimes_{R_{t(\sigma\op)}} M_{t(\sigma\op)})} {\bbQ / \bbZ}\\
& \cong \lim_{\sigma\op \in \calP_j(j,\bullet)} \Hom[\bbZ]{R_j \otimes_{R_{t(\sigma\op)}} M_{t(\sigma\op)}} {\bbQ / \bbZ} \\
& \cong \lim_{\sigma\op \in \calP_j(j,\bullet)} \Hom[R_{t(\sigma\op)}]{R_j}{\Hom[\bbZ]{M_{t(\sigma\op)}} {\bbQ / \bbZ}}\\
& = \lim_{\sigma\op \in \calP_j(j,\bullet)} \Hom[R_{t(\sigma\op)}]{R_j}{(M_{t(\sigma\op)})^+},
\end{align*}
and $({\varphi}_{j}^{M})^+$ is precisely ${\psi}_{j}^{M^+}$, the third statement is equivalent to that $M^+ \in \Psi(\sf Inj_{\bullet})$.
\end{prf*}

Then we have the following result, which generalizes a result by Enochs, Oyonarte and Torrecillas; see \cite[Theorem 3.7]{EOT04}.

\begin{thm} \label{flat right r module}
Suppose that $\calI$ is a direct category. Then there is an equality $$\F{\overline{\scrR} \lRep}=\Phi(\sf Flat_{\bullet}).$$
\end{thm}
\begin{prf*}
For the containment $\F{\overline{\scrR} \lRep} \subseteq \Phi(\sf Flat_{\bullet})$, we let $M$ be in $\F{\overline{\scrR}\lRep}$. Then one has $M^+\in\Psi(\sf Inj_{\bullet})$ by Lemma \ref{equivalent for flat} and \cite[Corollary 3.20]{DLLY}, and so $M$ is in $\Phi(\sf Flat_{\bullet})$ by Lemma \ref{corresponding relation}.

We then prove the containment $\Phi(\sf Flat_{\bullet})\subseteq\F{\overline{\scrR} \lRep}$. To this end, let $N$ be in $\Phi(\sf Flat_{\bullet})$. Then one has $N^+\in\Psi(\sf Inj_{\bullet})$ by Lemma \ref{corresponding relation}, and so $N^+$ is an injective object in $\widetilde{\scrR}\lRep$ by \cite[Corollary 3.20]{DLLY} as the index category $\calI\op$ for $\widetilde{\scrR}$ is inverse. This yields that $N$ is a flat object in $\overline{\scrR}\lRep$ by Lemma \ref{equivalent for flat}.
\end{prf*}

\section{Gorenstein flat model structure on $\overline{\scrR} \lRep$}
\label{Gorenstein flat}
\noindent
In this section, we construct the Gorenstein flat model structure on $\overline{\scrR} \lRep$, and give a characterization of Gorenstein flat objects in this category.

\begin{ipg}
For an arbitrary associative ring $A$, recall from Enochs, Jenda and Torrecillas \cite{EOT04} that a left $A$-module $M$ is called \emph{Gorenstein flat} if there is an exact sequence
\[
\mathbb{F}: \, \cdots \to F^{-1} \to F^0 \to F^{1} \to \cdots
\]
of flat left $A$-modules such that $M \cong \ker{(F^0 \to F^1)}$ and the sequence remains exact after applying the functor $\tp[A]{-}{E}$ for every injective right $A$-module $E$. Recently, \v{S}aroch and \v{S}t'ov\'{\i}\v{c}ek introduced \emph{projectively coresolved Gorenstein flat} left $A$-modules in \cite{SS20} by replacing flat modules in the above exact sequence with projective modules. Similarly one can define Gorenstein flat objects and projectively coresolved Gorenstein flat objects in $\overline{\scrR}\lRep$ using the categorical tensor product functors introduced in \ref{categorical tensor product} and injective objects in $\widetilde{\scrR}\lRep$.
\end{ipg}

\begin{nota}
Throughout this section, denote by $\scrR$ a representation of $\calI$ on {\sf Ring}  with $\scrR_i=R_i$ an associative ring for each $i\in\calI$ and $\scrR_\alpha$ a ring homomorphism for all $\alpha \in \Mor(\calI)$, and denote by
\begin{itemize}
\item $\GF{R_i \lMod}$ the subcategory of Gorenstein flat left $R_i$-modules;

\item $\PGF{R_i \lMod}$ the subcategory of projectively coresolved Gorenstein flat left $R_i$-modules;

\item $\sf GF_{\bullet}$ the family $\{ \GF{R_i \lMod} \}_{i \in \Ob(\calI)}$ of subcategories of $R_i \lMod$;

\item $\sf PGF_{\bullet}$ the family $\{\PGF{R_i \lMod}\}_{i \in \Ob(\calI)}$ of subcategories of $R_i \lMod$;

\item $\sf PGF ^\perp_{\bullet}$ the family $\{\PGF{R_i \lMod}^\perp\}_{i \in \Ob(\calI)}$ of subcategories of $R_i \lMod$;

\item $\overline{\scrR}$ the $\calI$-diagram of left module categories induced by $\scrR$;

\item $\GF{\overline{\scrR}\lRep}$ the subcategory of Gorenstein flat objects in $\overline{\scrR}\lRep$;

\item $\PGF{\overline{\scrR}\lRep}$ the subcategory of projectively coresolved Gorenstein flat objects in $\overline{\scrR}\lRep$.
\end{itemize}

Recall that a left $R_i$-module $M_i$ is called \emph{cotorsion} if $\Ext[R_i] {1}{F_i}{M_i} = 0$ for any flat left $R_i$-module $F_i$. Similarly, one can define cotorsion objects in $\overline{\scrR} \lRep$.
Denote by
\begin{itemize}
\item $\Cot{R_i \lMod}$ the subcategory of cotorsion left $R_i$-modules;

\item $\sf Cot_{\bullet}$ the family $\{\Cot{R_i \lMod}\}_{i \in \Ob(\calI)}$ of subcategories of $R_i \lMod$;

\item $\Cot{\overline{\scrR} \lRep}$ the subcategory of $\overline{\scrR} \lRep$ consisting of cotorsion objects.
\end{itemize}
\end{nota}
Recall that $\sf Proj_{\bullet}$ and $\sf Flat_{\bullet}$ denote the families $\{\Prj{R_i \lMod}\}_{i \in \Ob(\calI)}$ and $\{ \F {R_i \lMod} \}_{i \in \Ob(\calI)}$, respectively.

\begin{lem} \label{GF and PGF compatible}
The families $\sf Proj_{\bullet}$ and $\sf Flat_{\bullet}$ are compatible with respect to $\overline{\scrR}$. If $\scrR$ is flat, then the families $\sf GF_{\bullet}$ and $\sf PGF_{\bullet}$ are also compatible with respect to $\overline{\scrR}$.
\end{lem}

\begin{prf*}
The first statement is clear. For the second statement, we only deal with $\sf GF_{\bullet}$ as the argument also works for $\sf PGF_{\bullet}$ with small modifications. Since $\scrR$ is flat, $R_j$ is flat in $\rMod R_i$ for any $\alpha: i \to j \in \Mor(\calI)$. It is easy to check that $R_j \otimes_{R_i} G_i$ is Gorenstein flat in $R_j \lMod$ for any Gorenstein flat object $G_i \in R_i \lMod$; see \cite[Ascent table II(a)]{LWCHHl09}. Thus the family $\sf GF_{\bullet}$ is compatible with respect to $\overline{\scrR}$.
\end{prf*}

For each $i \in \Ob(\calI)$, by \v{S}aroch and \v{S}t'ov\'{\i}\v{c}ek \cite[Page 27]{SS20}, both
\begin{center}
$(\GF{R_i \lMod}, \, \PGF{R_i \lMod}^\perp, \, \Cot{R_i \lMod})$  and $(\PGF{R_i \lMod}, \, \PGF{R_i \lMod}^\perp, \, R_i \lMod)$
\end{center}
are hereditary Hovey triples in $R_i \lMod$ with
\begin{itemize}
\item[-] $\GF{R_i \lMod} \cap \PGF{R_i \lMod}^\perp = \F {R_i \lMod}$

\item[-] $\PGF{R_i \lMod}^\perp \cap \Cot{R_i \lMod} = \GF{R_i \lMod}^\perp$ and

\item[-] $\PGF{R_i \lMod} \cap \PGF{R_i \lMod}^\perp = \Prj{R_i \lMod}$.
\end{itemize}

\begin{thm}\label{ht induce ht 12}
Suppose that $\calI$ is a left rooted quiver and $\scrR$ is flat. Then both
\[
(\Phi(\sf GF_{\bullet}), \, \overline{\scrR}\lRep_{\sf PGF^\perp_{\bullet}}, \, \overline{\scrR} \lRep_{\sf Cot_{\bullet}}) \text{ and } (\Phi(\sf PGF_{\bullet}), \, \overline{\scrR} \lRep_{\sf PGF^\perp_{\bullet}}, \, \overline{\scrR}\lRep)
\]
are hereditary Hovey triples in $\overline{\scrR} \lRep$.
\end{thm}
\begin{prf*}
Since $\scrR$ is flat, one gets that the $\calI$-diagram $\overline{\scrR}$ is exact. Thus by Theorem \ref{ht induce ht 1} and Lemma \ref{GF and PGF compatible}, both $(\Phi(\sf GF_{\bullet}), \, \overline{\scrR}\lRep_{\sf PGF^\perp_{\bullet}}, \, \overline{\scrR} \lRep_{\sf Cot_{\bullet}})$ and $(\Phi(\sf PGF_{\bullet}), \, \overline{\scrR} \lRep_{\sf PGF^\perp_{\bullet}}, \, \overline{\scrR}\lRep)$ are hereditary Hovey triples.
\end{prf*}

\begin{cor} \label{some cps over direct}
Suppose that $\calI$ is a left rooted quiver, and $\scrR$ is flat. Then
\begin{center}
$(\Phi(\sf Flat_{\bullet}), \, \overline{\scrR} \lRep_{\sf Cot_{\bullet}})$, $(\Phi(\sf GF_{\bullet}), \, \overline{\scrR} \lRep_{\sf GF^\perp_{\bullet}})$ and $(\Phi(\sf PGF_{\bullet}), \, \overline{\scrR} \lRep_{\sf PGF^{\perp}_{\bullet}})$
\end{center}
are complete and hereditary cotorsion pairs in $\overline{\scrR} \lRep$.
\end{cor}

\begin{prf*}
Note that both $\GF{R_i \lMod}$ and $\F {R_i \lMod}$ are closed under small colimits for each $i \in \Ob(\calI)$. Then by Proposition \ref{com and her cp for q tilte} and Lemma \ref{GF and PGF compatible}, we conclude that $\Phi(\sf GF_{\bullet})\cap \overline{\scrR}\lRep_{\sf PGF^\perp_{\bullet}}=\Phi(\sf Flat_{\bullet})$, and so the first pair is a complete and hereditary cotorsion pair in $\overline{\scrR}\lRep$. Theorem \ref{ht induce ht 12} tells us that the second and the third pair are also complete and hereditary cotorsion pairs in $\overline{\scrR}\lRep$.
\end{prf*}

Under some conditions the subcategories $\Cot{\overline{\scrR} \lRep}$ and $\overline{\scrR} \lRep_{\sf Cot_{\bullet}}$ coincide.

\begin{prp} \label{cotorsion object}
Suppose that $\calI$ is a left rooted quiver, and $\scrR$ is flat. Then there is an equality $\Cot{\overline{\scrR} \lRep} = \overline{\scrR} \lRep_{\sf Cot_{\bullet}}$.
\end{prp}

\begin{prf*}
Note that the notion of flat objects given in Definition \ref{definition of flat} is indeed the categorical flat objects in $\overline{\scrR} \lRep$; see \cite{Oberst70}. It follows that $(\F{{\overline{\scrR} \lRep}}, \, \Cot{{\overline{\scrR} \lRep}})$ is a cotorsion pair in $\overline{\scrR} \lRep$. On the other hand, by Corollary \ref{some cps over direct}, $(\Phi(\sf Flat_{\bullet}), \, \overline{\scrR} \lRep_{\sf Cot_{\bullet}})$
is also a cotorsion pair in $\overline{\scrR} \lRep$. However, Theorem \ref{flat right r module} tells us that $\F{\overline{\scrR} \lRep} = \Phi(\sf Flat_{\bullet})$. Thus $\Cot{\overline{\scrR} \lRep} = \overline{\scrR} \lRep_{\sf Cot_{\bullet}}$.
\end{prf*}

Let $\calY$ be a subcategory of an abelian category $\calA$. Recall that a subcategory $\calV$ is called a \emph{cogenerator} \cite{MAsROB89} for $\calY$ if for any object $Y \in \calY$,
there exists a short exact sequence $0 \to  Y \to V \to Y' \to 0$ with $V \in \calV$ and
$Y' \in \calY$. The next result will be applied in the proof of
Theorem \ref{Gorenstein flat r module}.

\begin{lem} \label{flat cogenertor}
Suppose that $\calI$ is a left rooted quiver and $\scrR$ is flat. Then the following hold.
\begin{prt}
\item $\Phi(\sf Flat_{\bullet})$ is a cogenerator for $\Phi(\sf GF_{\bullet})$;

\item $\Phi(\sf Proj_{\bullet})$ is a cogenerator for $\Phi(\sf PGF_{\bullet})$.
\end{prt}
\end{lem}

\begin{prf*}
We only prove statement (a) since the other one can be proved similarly. Let $M$ be an object in $\Phi(\sf GF_{\bullet})$. Since $(\Phi(\sf GF_{\bullet}), \, \overline{\scrR} \lRep_{\sf GF^\perp_{\bullet}})$ is a complete cotorsion pair in $\overline{\scrR} \lRep$ by Corollary \ref{some cps over direct}, there exists a short exact sequence
$0 \to M \to F \to M' \to 0$ with  $F \in \Phi(\sf GF_{\bullet}) \cap \overline{\scrR} \lRep_{\sf GF^\perp_{\bullet}} $ and  $M' \in \Phi(\sf GF_{\bullet})$
in $\overline{\scrR} \lRep$. Thus it is enough to show that $F \in \Phi(\sf Flat_{\bullet})$. Indeed, since
\[
(\Phi(\sf GF_{\bullet}), \, \overline{\scrR} \lRep_{\sf PGF^\perp_{\bullet}}, \, \overline{\scrR} \lRep_{\sf Cot_{\bullet}})
\]
is a Hovey triple in $\overline{\scrR} \lRep$ by Proposition \ref{ht induce ht GI 1}, we have
\[
\Phi(\sf GF_{\bullet}) \cap \overline{\scrR} \lRep_{\sf PGF^\perp_{\bullet}} = \F{\overline{\scrR} \lRep}
\]
as $(\F{{\overline{\scrR} \lRep}}, \, \Cot{\overline{\scrR} \lRep})$ is a cotorsion pair in $\overline{\scrR} \lRep$ and $\Cot{\overline{\scrR} \lRep} = \overline{\scrR} \lRep_{\sf Cot_{\bullet}}$ by Proposition \ref{cotorsion object}. By Theorem \ref{flat right r module}, $\F{\overline{\scrR} \lRep} = \Phi(\sf Flat_{\bullet})$, so
\[
\Phi(\sf GF_{\bullet}) \cap \overline{\scrR} \lRep_{\sf GF^\perp_{\bullet}} = \Phi(\sf Flat_{\bullet})
\]
and $F \in \Phi(\sf Flat_{\bullet})$.
\end{prf*}

In the following lemma, we collect some elementary properties of Gorenstein flat objects in $\overline{\scrR} \lRep$, which are quite similar to those of Gorenstein injective objects in $\widetilde{\scrR} \lRep$; see Lemma \ref{on side for g injectives}.

\begin{lem} \label{on side for g flat}
Suppose that $\calI$ is a left rooted quiver and $\scrR$ is flat. Let $M$ be a Gorenstein flat object in $\overline{\scrR} \lRep$. Then for each $i \in \Ob(\calI)$, ${\varphi}_i^{M}: \coprod_{\theta \in \calI(\bullet,i)} (R_i \otimes_{R_{s(\theta)}} M_{s(\theta)}) \to M_i$ is a monomorphism with $\coker({\varphi}_i^{M})$ Gorenstein flat in $R_i \lMod$.
That is, there is a containment $\GF{\overline{\scrR} \lRep} \subseteq \Phi(\sf GF_{\bullet})$.
\end{lem}

\begin{prf*}
Fix $i \in \Ob(\calI)$. We define a functor $\C_i: \overline{\scrR}\lRep\to R_i\lMod$ sending a representation $M\in\overline{\scrR}\lRep$ to $\coker({\varphi}_i^{M})$; see \cite[Corollary 2.15]{DLLY}. Since $M$ is a Gorenstein flat object in $\overline{\scrR}\lRep$, there exists an exact sequence
\[
\mathbb{F}: \, \cdots \to F^{-1} \to F^0 \to F^1 \to \cdots
\]
of flat objects in $\overline{\scrR}\lRep$ such that $M \cong \ker{(F^0 \to F^1)}$ and the sequence $\mathbb{F}$ remains exact after applying the functor $E\otimes_R-$ for every injective object $E \in \widetilde{\scrR} \lRep$. For all integers $j$, since $F^j \in \Phi(\sf Flat_{\bullet})$ by Theorem \ref{flat right r module}, there is a short exact sequence
\[
0 \to \coprod_{\theta \in \calI(\bullet,i)} (R_i \otimes_{R_{s(\theta)}} F^j_{s(\theta)}) \to F^j_i \to \C_i(F^j) \to 0
\]
in $R_i \lMod$ with $\C_i(F^j)$ flat. For any arrow $\theta \in \calI(\bullet,i)$, we have an exact sequence
\[
\mathbb{F}_{s(\theta)}: \, \cdots \to F^{-1}_{s(\theta)} \to F^0_{s(\theta)} \to F^{1}_{s(\theta)} \to \cdots
\]
in $\overline{\scrR}_{s(\theta)}$. Since $R_i$ is a flat right $R_{s(\theta)}$-module as $R$ is flat by assumption, the sequence $R_i \otimes_{R_{s(\theta)}} \mathbb{F}_{s(\theta)}$ is exact, so the sequence
$\coprod_{\theta \in \calI(\bullet,i)} (R_i \otimes_{R_{s(\theta)}} \mathbb{F}_{s(\theta)})$ is also exact. Consequently, we obtain the commutative diagram
\[
\xymatrix@C=20pt@R=15pt{
 & \vdots \ar[d] & \vdots \ar[d] & \vdots \ar[d] \\
0 \ar[r] & \coprod_{\theta \in \calI(\bullet,i)} (R_i \otimes_{R_{s(\theta)}} F^{-1}_{s(\theta)}) \ar[d] \ar[r] & F^{-1}_i \ar[d] \ar[r] & \C_i(F^{-1}) \ar[d] \ar[r] & 0 \\
0 \ar[r] & \coprod_{\theta \in \calI(\bullet,i)} (R_i \otimes_{R_{s(\theta)}} F^0_{s(\theta)}) \ar[d] \ar[r] & F^0_i \ar[d] \ar[r] & \C_i(F^0) \ar[d] \ar[r] & 0 \\
0 \ar[r] & \coprod_{\theta \in \calI(\bullet,i)} (R_i \otimes_{R_{s(\theta)}} F^{1}_{s(\theta)}) \ar[d] \ar[r] & F^{1}_i \ar[d] \ar[r] & \C_i(F^{1}) \ar[d] \ar[r] & 0 \\
 & \vdots & \vdots & \vdots
}
\]
with exact rows and columns, which induces the short exact sequence
\[
0 \to \coprod_{\theta \in \calI(\bullet,i)} (R_i \otimes_{R_{s(\theta)}} {M_{s(\theta)}}) \overset{\varphi_i^M} \longrightarrow M_i \to \C_i(M) \to 0.
\]
Therefore, to complete the proof, it remains to show that $\C_i(M)$ is Gorenstein flat in $R_i \lMod$, that is, the sequence $I_i \otimes_{R_i} \C_i(\mathbb{F})$ is exact for each injective right $R_i$-module $I_i$.

For any injective object $E \in \widetilde{\scrR} \lRep$, the sequence $\Hom[{\widetilde{\scrR} \lRep}] {E} {\mathbb{F}^+}$ is exact as
\[
\Hom[{\widetilde{\scrR} \lRep}]{E}{\mathbb{F}^+} \cong (E\otimes_R \mathbb{F})^+
\]
by Lemma \ref{adjoint-iso for tensor product}. By an argument similar to the one used in the proof of Lemma \ref{on side for g injectives}, we conclude that the sequence $\Hom[R_i] {I_i} {\K_i(\mathbb{F}^+)}$ is exact. Now the conclusion follows from the above isomorphism and the observation that $\K_i(\mathbb{F}^+)$ is exactly $\C_i(\mathbb{F})^+$.
\end{prf*}

According to \thmcite[3.6]{HHl04a}, if $M_i$ is Gorenstein flat in $R_i \lMod$, then $M_i^+$ is Gorenstein injective in $\rMod R_i$. The converse statement is also true whenever $R_i$ is right coherent. Relying on this fact, one can prove the following result by an argument similar to the one used in the proof of Lemma \ref{corresponding relation}.

\begin{lem} \label{corresponding relation for gflat}
Suppose that $\calI$ is a left rooted quiver, and let $M$ be an object in $\overline{\scrR} \lRep$. If $M$ is contained in $\Phi(\sf GF_{\bullet})$, then $M^+$ is contained in the subcategory $\Psi(\sf GI_{\bullet})$ of $\widetilde{\scrR} \lRep$. Furthermore, the converse statement holds if $R_i$ is right coherent for every $i \in \Ob(\calI)$.
\end{lem}

Now we can give characterizations of Gorenstein flat objects and projectively coresolved Gorenstein flat objects in $\overline{\scrR} \lRep$.

\begin{thm} \label{Gorenstein flat r module}
Suppose that $\calI$ is a left rooted quiver, and $\scrR$ is flat. Then there are equalities
$$\GF{\overline{\scrR} \lRep} = \Phi({\sf GF}_{\bullet})\ \  \text{and}\ \ \PGF{\overline{\scrR} \lRep} = \Phi(\sf PGF_{\bullet}).$$
\end{thm}

\begin{prf*}
We only show the first equality since the second one can be proved similarly. Lemma \ref{on side for g flat} tells us that $\GF{\overline{\scrR} \lRep} \subseteq \Phi(\sf GF_{\bullet})$, so we only need to show the other inclusion.

Take $M \in \Phi(\sf GF_{\bullet})$. By Theorem \ref{flat right r module} and Lemma \ref{flat cogenertor}(a), $\F{\overline{\scrR} \lRep} = \Phi(\sf Flat_{\bullet})$ is a cogenerator for $\Phi(\sf GF_{\bullet})$, so there is a short exact sequence $0 \to M \to F^0 \to M^{1} \to 0$ in $\overline{\scrR} \lRep$ with $F^0 \in \F{\overline{\scrR} \lRep}$ and $M^{1} \in \Phi(\sf GF_{\bullet})$, which induces a short exact sequence
\[
0 \to (M^{1})^+ \to (F^0)^+ \to M^+ \to 0
\]
in $\widetilde{\scrR} \lRep$. By Lemma \ref{corresponding relation for gflat}, $(M^{1})^+$ is contained in $\Psi(\sf GI_{\bullet})$. Therefore, by Theorem \ref{Gorenstein in r module}, $(M^{1})^+$ is Gorenstein injective in $\widetilde{\scrR} \lRep$, and hence the sequence
\[
0 \to \Hom[{\widetilde{\scrR} \lRep}]{E} {(M^{1})^+} \to \Hom[{\widetilde{\scrR} \lRep}]{E}{(F^0)^+} \to \Hom[{\widetilde{\scrR} \lRep}]{E}{M^+} \to 0
\]
is exact for any injective object $E \in \widetilde{\scrR} \lRep$. It follows from Lemma \ref{adjoint-iso for tensor product} that the sequence
\[
0 \to E \otimes_R M \to E \otimes_R F^0 \to E \otimes_R M^{1} \to 0
\]
is also exact. Replacing $M$ by $M^{1}$, recursively one gets an exact sequence
\[
0 \to M \to F^0 \to F^{1} \to \cdots \quad \quad \quad{(\dag)}
\]
in $\overline{\scrR} \lRep$ with $F^j \in \F{\overline{\scrR} \lRep}$ for all integers $j \geqslant 0$ such that the sequence $(\dag)$  remains exact after applying the functor $E \otimes_R -$.

On the other hand, note that $\overline{\scrR} \lRep$ has enough projectives by Theorem \ref{Rep is grothendieck}, so there exists an exact sequence
\[
\cdots \to F^{-2} \to F^{-1} \to M \to 0\quad \quad \quad (\ddag)
\]
in $\overline{\scrR} \lRep$ with $F^j \in \Prj{\overline{\scrR} \lRep} \subseteq \F{\overline{\scrR} \lRep}$ for all integers $j < 0$. Since $\Phi(\sf GF_{\bullet})$ is closed under taking kernels of epimorphisms by Corollary \ref{some cps over direct}, all kernels of the sequence $(\ddagger)$ belong to $\Phi(\sf GF_{\bullet})$ by noting that each $F^j \in \Phi(\sf Flat_{\bullet})$ and $\Phi(\sf Flat_{\bullet}) \subseteq \Phi(\sf GF_{\bullet})$ clearly. Using a similar argument as before we conclude that the sequence $(\ddag)$ remains exact after applying the functor $E \otimes_R -$.

Assembling the exact sequences $(\dagger)$ and $(\ddagger)$, we deduce that $M$ is in $\GF{\overline{\scrR} \lRep}$, as desired.
\end{prf*}

An immediate consequence of Corollary \ref{some cps over direct} and Theorem \ref{Gorenstein flat r module} is:

\begin{cor} \label{cp for pgf}
Suppose that $\calI$ is a left rooted quiver, and $\scrR$ is flat. Then
\begin{center}
$(\GF{\overline{\scrR} \lRep}, \, \overline{\scrR} \lRep_{\sf GF^{\perp}_{\bullet}})$ and $(\PGF{\overline{\scrR} \lRep}, \, \overline{\scrR} \lRep_{\sf PGF^{\perp}_{\bullet}})$
\end{center}
are complete and hereditary cotorsion pairs in $\overline{\scrR} \lRep$.
\end{cor}

The following result provides a Gorenstein flat and a projectively coresolved Gorenstein flat model structures on $\overline{\scrR} \lRep$.

\begin{cor} \label{gflat model structor}
Suppose that $\calI$ is a left rooted quiver, and $\scrR$ is flat. Then there exist hereditary Hovey triples
\[
(\GF{\overline{\scrR} \lRep}, \, \PGF{\overline{\scrR} \lRep}^\perp, \, \Cot{\overline{\scrR} \lRep}) \text{ and } (\PGF{\overline{\scrR} \lRep},\, \PGF{\overline{\scrR} \lRep}^\perp,\, \overline{\scrR}\lRep)
\]
in $\overline{\scrR} \lRep$.
\end{cor}

\begin{prf*}
By Theorem \ref{ht induce ht 12}, the triples
\[
(\Phi({\sf GF}),\, \overline{\scrR} \lRep_{\sf PGF^\perp_{\bullet}},\, \overline{\scrR} \lRep_{\sf Cot_{\bullet}}) \text{ and } (\Phi(\sf PGF_{\bullet}),\, \overline{\scrR} \lRep_{\sf PGF^\perp_{\bullet}}, \, \overline{\scrR}\lRep)
\]
are hereditary Hovey triples in $\overline{\scrR} \lRep$. Note that $\Phi(\sf GF_{\bullet}) = \GF{\overline{\scrR} \lRep}$ and $\Phi(\sf PGF_{\bullet}) = \PGF{\overline{\scrR} \lRep}$ by Theorem \ref{Gorenstein flat r module}, $\overline{\scrR} \lRep_{\sf PGF^\perp_{\bullet}} = \PGF{\overline{\scrR} \lRep}^\perp$ by Corollary \ref{cp for pgf}, and $\overline{\scrR} \lRep_{\sf Cot_{\bullet}} = \Cot{\overline{\scrR} \lRep}$ by Proposition \ref{cotorsion object}. Consequently, the conclusion follows.
\end{prf*}

\bigskip
\section*{Acknowledgments}%%%%%%%%%%%%%%%%%%%%%%%%%%%%%%%%%%%%%%%%%%%%
\noindent
We thank Fei Xu for valuable discussions related to this work. Z.X. Di was partially supported by the National Natural Science Foundation of China (Grant No. 12471034), the Scientific Research Funds of Fujian Province (Grant No. 605-52525002), and the Scientific Research Funds of Huaqiao University (Grant No. 605-50Y22050). L.P. Li was partially supported by the National Natural Science Foundation of China (Grant No. 12171146). L. Liang was partially supported by the National Natural Science Foundation of China (Grant No. 12271230) and the Foundation for Innovative Fundamental Research Group Project of Gansu Province (Grant No. 25JRRA805). N.N. Yu was partially supported by the Natural Science Foundation of Xiamen Municipality (Grant No. 3502Z202473007), the Natural Science Foundation of Fujian Province (Grant No. 2024J01027), and the National Natural Science Foundation of China (Grant Nos. 12131018 and 12471121).

\bibliographystyle{amsplain-nodash}

%\bibliographystyle{amsplain}
%\bibliography{./+references}

\def\cprime{$'$}
  \providecommand{\arxiv}[2][AC]{\mbox{\href{http://arxiv.org/abs/#2}{\sf
  arXiv:#2 [math.#1]}}}
  \providecommand{\oldarxiv}[2][AC]{\mbox{\href{http://arxiv.org/abs/math/#2}{\sf
  arXiv:math/#2
  [math.#1]}}}\providecommand{\MR}[1]{\mbox{\href{http://www.ams.org/mathscinet-getitem?mr=#1}{#1}}}
  \renewcommand{\MR}[1]{\mbox{\href{http://www.ams.org/mathscinet-getitem?mr=#1}{#1}}}
\providecommand{\bysame}{\leavevmode\hbox to3em{\hrulefill}\thinspace}
\providecommand{\MR}{\relax\ifhmode\unskip\space\fi MR }
% \MRhref is called by the amsart/book/proc definition of \MR.
\providecommand{\MRhref}[2]{%
  \href{http://www.ams.org/mathscinet-getitem?mr=#1}{#2}
}
\providecommand{\href}[2]{#2}

\end{document}